\documentclass[reqno]{amsart}
\usepackage{latexsym}
\theoremstyle{plain}
\usepackage{txfonts}
\usepackage[sc]{mathpazo}
\def\Xint#1{\mathchoice
   {\XXint\displaystyle\textstyle{#1}}%
   {\XXint\textstyle\scriptstyle{#1}}%
   {\XXint\scriptstyle\scriptscriptstyle{#1}}%
   {\XXint\scriptscriptstyle\scriptscriptstyle{#1}}%
   \!\int}

\def\XXint#1#2#3{{\setbox0=\hbox{$#1{#2#3}{\int}$}
     \vcenter{\hbox{$#2#3$}}\kern-.5\wd0}}

\newcommand{\D}{\mathcal{D}}

\newcommand{\chara}{1\!\!1}

\newcommand{\cha}{1\!\!1}

\newcommand{\na}{D}

\newcommand{\ui}{\tilde{u}}

\newcommand{\lt}{\left}
\newcommand{\rt}{\right}
\newcommand{\nl}{\newline}

\newcommand{\nn}{\nonumber}

\newcommand{\lm}{\lambda}

\newcommand{\qd}{\quad}

\newcommand{\ep}{\epsilon}

\newcommand{\II}{\mathcal{I}}

\newcommand{\PPI}{\mathcal{P}}

\newcommand{\GI}{\mathcal{I}}

\newcommand{\SI}{\mathcal{S}}

\newcommand{\CI}{\mathcal{C}}

\newcommand{\OI}{\mathcal{O}}
\newcommand{\BI}{\mathcal{B}}

\newcommand{\ti}{\tilde}

\newcommand{\R}{\mathrm {I\!R}}

\newcommand{\dia}{\diamondsuit}

\newcommand{\mucsta}{\mu^{38Q}}

\newcommand{\mucstc}{9C_p (1352)^{18Q} Q^{40Q}}
\newcommand{\mucon}{(2\times 10^{10}(Q+1)^6)^{-6(Q+1)}}
\newcommand{\cfc}{3.7\times 10^{62}}
\newcommand{\gammacon}{\lt(\frac{\mu}{2000 (Q+1)^2}\rt)^{6(Q+1)}}
\newtheorem{a1}{Lemma}
\newtheorem{a2}{Theorem}
\newtheorem{a3}{Conjecture}

\newtheorem{a5}{Proposition}

\newtheorem{deff}{Definition}
\newtheorem{quest}{Question}
\newtheorem*{thm}{Theorem}

\theoremstyle{remark}

\begin{document}
\title[Rigidity of pairs quasiregular mappings whose symmetric part of gradient are close]
{Rigidity of pairs of quasiregular mappings whose symmetric part of gradient are close}
\author{Andrew Lorent}
\address{Mathematics Department\\University of Cincinnati\\2600 Clifton Ave.\\ Cincinnati OH 45221 }
\email{lorentaw@uc.edu}
\maketitle

\begin{abstract} For $A\in M^{2\times 2}$ let $S(A)=\sqrt{A^T A}$, i.e.\ the symmetric part of the 
polar decomposition of $A$. 
We consider the relation between two quasiregular mappings whose symmetric part of gradient are close. Our main 
result is the following. Suppose $v,u\in W^{1,2}(B_1(0):\R^2)$ 
are $Q$-quasiregular mappings with $\int_{B_1(0)} \det(Du)^{-p} dz\leq C_p$ for 
some $p\in (0,1)$ and 
$\int_{B_1(0)} \lt|Du\rt|^2 dz\leq 1$. There exists constant $M>1$ such that if  
$\int_{B_1(0)} \lt|S(Du)-S(Dv)\rt|^2 dz=\ep$ then 
$$
\int_{B_{\frac{1}{2}}(0)} \lt|Dv-R Du\rt| dz\leq c\CI_p^{\frac{1}{p}}\ep^{\frac{p^3}{M Q^5\log(10 C_p Q)}}\text{ for some }R\in SO(2).
$$
Taking $u=Id$ we obtain a special case of the quantitative rigidity 
result of Friesecke, James and M\"{u}ller \cite{fmul}. Our main result can be considered as a first step in a new line of generalization of Theorem 1 of \cite{fmul} in which $Id$ is replaced by a mapping of non-trivial degree.
\end{abstract}

Rigidity and stability of differential inclusions is a classical subject. Reshetnyak's monograph \cite{res1} is devoted to proving a 
quantitative stability result generalizing Liouville's classic theorem \cite{lou} that solutions of the differential inclusion $Du\in CO_{+}(n):=\lt\{\lm R:\lm>0,R\in SO(n)\rt\}$, $n\geq 3$ are affine or Mobius. Korn's inequality is an optimal quantitative stability result for the fact that the differential inclusion 
$Du\in \mathrm{Skew}(n\times n):=\lt\{M\in M^{n\times n}:M^T=-M\rt\}$ is satisfied only by an affine map. 

This subject has received considerable impetus from the work of Friesecke, James and M\"{u}ller \cite{fmul} who proved an optimal quantitative stability 
result for the corollary to Liouville's theorem that states solutions to the differential inclusion $Du\in SO(n)$ are affine.  
\begin{a2}[Friesecke, James and M\"{u}ller, 2002]
\label{MT} For every bounded open connected Lipschitz domain $U\subset \R^n$, $n\ge 2$, 
and every $q>1$, there exists a constant $C=C(U,q,n)$
such that writing $K := SO(n)$, 
$$
\inf_{R\in K} \label{mthm1} \| D v-R\|_{L^q\lt(U\rt)}\leq C\|d\lt(D v,K\rt)\|_{L^q\lt(U\rt)}
\quad\mbox{ for every } v\in W^{1,q}(U;\R^n).
$$
\end{a2}
Previously strong partial results controlling the function 
(rather than the gradient) have been established by John \cite{john}, Kohn \cite{kohn}. 

The simplicity of the statement of Theorem \ref{MT} can lead to the strength of the advance that is represented by this theorem being overlooked. It is rare in contemporary 
research in analysis to prove a new and deep result about elementary 
mathematical objects; Theorem \ref{MT} is exactly such a result. It has 
had wide application in applied analysis and is one of the main tools used to 
make a rigorous and complete analysis of the multiple thin shell theories in classical 
elasticity \cite{fmul}, \cite{fmul1}, \cite{fmul2}. Beyond this it has the merit of being a statement whose significance would be clear to mathematicians of two hundred years ago.

A number of works have extended Theorem \ref{MT} to cover various larger classes of matrices than $SO\lt(n\rt)$. Faraco and Zhong proved 
the corresponding result with
$K=\Pi SO\lt(n\rt)$ where $\Pi \subset \R_{+}\backslash \lt\{0\rt\}$ is a compact set, 
\cite{fac}. 
Chaudhuri and M\"uller \cite{chau} and later Delellis and  Szekelyhidi \cite{las1} considered a set of the form $K=SO\lt(n\rt)A\cup SO\lt(n\rt)B$ where $A$ and $B$ are \em strongly incompatible \rm in the sense of Matos \cite{matos}.

In this paper, following an approach started by Ciarlet and Mardare \cite{cia} and also suggested by M\"{u}ller, we start a different line of generalization of Theorem \ref{MT}. The initial 
observation is that Theorem \ref{MT} is a special case of the following question. Recall we defined $S(M)=\sqrt{M^T M}$ to be the symmetric part of a matrix. 

\begin{quest}
\label{q1}
If $\Omega\subset \R^n$ is a connected domain and $u,v\in W^{1,2}(\Omega)$, $\mathrm{det}(Du)>0$, $\mathrm{det}(Dv)>0$ and 
$\int_{\Omega} \lt|S(Du)-S(Dv)\rt|^2 dx=\ep$ does this imply there exists $R\in SO(n)$ such that $\int_{\Omega} \lt|Du-R Dv\rt|^2 dx\leq \delta$ where $\delta$ is some small quantity 
depending on $\ep$. 
\end{quest}

It turns out that the answer to Question \ref{q1} is no, even in the "absolute" version of this question where $\ep=0$, see Example 1 \cite{lor14} or see the example in Section 4, \cite{cia}. 
For a positive result for the case where 
$\ep=0$ it suffices to consider the class of functions of integrable dilatation as shown in Theorem 1 \cite{lor14} (or see Theorem 1 of \cite{lor27} for a more general result). Theorem 1 of \cite{lor27} and the 2d version of Theorem 1 of \cite{lor14} are sharp in the sense that no result of this kind is possible outside the space of mappings of integrable dilatation.

In this paper we will provide a positive answer to Question \ref{q1} for pairs of Quasiregular mappings in two dimensions. Note in Theorem \ref{T1} and throughout the paper a ball of radius $r$ centred on zero will be denoted $B_r$.
\begin{a2}
\label{T1}
Suppose $v,u\in W^{1,2}(B_1:\R^2)$ 
are $Q$-quasiregular mappings with $\int_{B_1} \det(Du)^{-p} dz\leq C_p$ for 
some $p\in (0,1)$ and 
$\int_{B_1} \lt|Du\rt|^2 dz\leq 1$. If   
\begin{equation}
\label{fg1}
\int_{B_1} \lt|S(Du)-S(Dv)\rt|^2 dz=\ep
\end{equation}
then there exists $R\in SO(2)$ such that 
\begin{equation}
\label{fest1}
\int_{B_{\frac{1}{2}}} \lt|Dv-R Du\rt| dz\leq c C_p^{\frac{1}{p}}\ep^{\frac{p^3}{10^{8} Q^5\log(10 C_p Q)}}. 
\end{equation}
\end{a2}
Theorem \ref{T1} to a certain extent shares the property that Theorem \ref{MT} has of being a new and interesting statement about the classical objects of mathematical analysis. The credit for this however is largely due to Theorem 
\ref{MT} as the methods of proof of this theorem are used in an essential 
way in the proof of Theorem \ref{T1}. In this author's opinion there are a number of results in the area of classical Quasiconformal analysis that can be harvested 
by use of the ideas in the proof of Theorem \ref{MT}, Theorem \ref{T1} is just 
one of them. Note if we take $u=Id$ hypothesis (\ref{fg1}) is exactly $\int_{B_1} d^2(Dv,SO(2)) dz=\ep$ and the 
conclusion is $\int_{B_{\frac{1}{2}}} \lt|Dv-R\rt| dz\leq c\ep^{\frac{p^3}{10^{8} \log(10 )}}$ for 
some $R\in SO(2)$. While this is much weaker than Theorem \ref{MT} it is still a result that was not known 
prior to the publication of \cite{fmul}. In some sense the line of generalization that this paper 
contributes to is the desire to replace $Id$ by a mapping of non-trivial degree.

Ciarlet and Mardare were motivated to study Question \ref{q1} as part of a program to develop a theory 
of elasticity based on study the "Cauchy Green" tensor $Du^T Du$ of a deformation $u$, \cite{cia}, \cite{cia1}, \cite{cia2}. They 
proved a version of Theorem \ref{T1} for $C^1$ mappings with the property that $\det(Du)>0$ everywhere in the 
domain and the constant $c$ in (\ref{fest1}) depends on $u$. Their method was again to apply Theorem \ref{MT}, this 
will be sketched in the next section.  

Theorem \ref{T1} is clearly suboptimal however we believe the power of $\ep$ in inequality (\ref{fest1}) is 
of the right form in the sense that the power decreases as the degree of the mapping $u$ increases or as $Q$ increases. As the dependence on the degree is a key issue an 
example showing the dependence will be presented in \cite{lor39}. We give a sketch of the construction of the 
example in Section \ref{example}.

\section{Proof Sketch}

\subsection{Absolute case with global invertibility.}
\label{skk1}

First suppose we have 
$C^1$ functions $u,v$ where $u$ is globally invertible and $S(Du)=S(Dv)$ everywhere. By polar decomposition we have 
$A=R(A)S(A)$ for some $R(A)\in SO(n)$. Form $w(z)=v(u^{-1}(z))$ and 
note that 
\begin{eqnarray}
\label{opz1}
\na w(x)&=&\na v(u^{-1}(x))(\na u(u^{-1}(x)))^{-1}\nn\\
&=&R(\na v(u^{-1}(x)))\lt(R(\na u(u^{-1}(x)))\rt)^{-1}\in SO(n)\nn
\end{eqnarray}
by Liouville's theorem it is clear there exists $R\in SO(n)$ such that 
$\na w(z)=R$ for all $z\in \Omega$. Thus 
\begin{equation}
\label{opx1}
\na v=R \na u\text{ on }\Omega 
\end{equation}
and result is established. 

\subsection{Quantitative case with global invertibility.}
\label{QG}

Now assume $u,v$ are $C^1$ and $u$ is globally invertible and $\int_{B_1} \lt|S(Du)-S(Dv)\rt|^2 dz=\ep$ and 
$\inf\lt\{\det(Du(z)):z\in \Omega\rt\}>0$. Apart from where $\lt|Du\rt|\sim 0$ and $\lt|Du\rt|\sim 0$ we know $\lt|(S(Du(z)))^{-1}-(S(Dv(z)))^{-1}\rt|\approx \lt|S(Du(z))-S(Dv(z))\rt|$ and hence letting 
$$
E(z)=\lt(S(Du(z))\rt)^{-1}-\lt(S(Dv(z))\rt)^{-1}
$$
we have 
\begin{eqnarray}
Dw(x)&=&
R(\na v(u^{-1}(x)))S(\na v(u^{-1}(x))\lt(S(\na u(u^{-1}(x))\rt)^{-1}\lt(R(\na u(u^{-1}(x))\rt)^{-1}\nn\\
&=&
R(\na v(u^{-1}(x)))S(\na v(u^{-1}(x))\lt(\lt(S(\na v(u^{-1}(x))\rt)^{-1}+E(u^{-1}(x)) \rt)\lt(R(\na u(u^{-1}(x))\rt)^{-1}\nn\\
&=&R(\na v(u^{-1}(x)))\lt(R(\na u(u^{-1}(x))\rt)^{-1}\nn\\
&~&\qd\qd+R(\na v(u^{-1}(x)))S(\na v(u^{-1}(x))E(u^{-1}(x)) \lt(R(\na u(u^{-1}(x))\rt)^{-1}.\nn
\end{eqnarray}
So 
\begin{eqnarray}
\int_{u(\Omega)} d^{2}\lt(Dw(z),SO(2)\rt) dz&\leq& c\lt(\int_{\Omega} \lt|Dv(z)\rt|^2 dz\rt)^{\frac{1}{2}}\lt(\int_{\Omega} E(z) dz\rt)^{\frac{1}{2}}\nn\\
&\leq& c\ep.
\end{eqnarray}
So applying Theorem \ref{MT} we have that there is constant $\CI=\CI(u)$ such that 
$$
\int_{u(\Omega)} \lt|Dw(z)-R_0\rt|^2 dz\leq \CI \ep
$$
and unwrapping gives the estimate we seek, however with a constant depending on $u$.

\subsection{Sketch of the General case}

Our problem is that we do not have global invertibility and we would like an estimate 
that depends on $u$ in a more explicit way. Under the 
hypothesis that the mappings $u,v$ are $Q$-quasiregular  we know that $u$ is locally invertible at all but countably many points, but we have no estimates of the size of the of neighbourhoods of invertibility. If we 
wanted to prove an estimate of the form (\ref{fest1}) where the constant $c$ depended on $u$ we could 
patch together neighbourhoods of invertibility so long as we knew the "size" of the neighbourhoods were bounded below on all 
compact subdomains. Under the hypothesis $\det(Du)>0$ everywhere for a $C^1$ function $u$ this is true 
and this is how Ciarlet and Mardare established their estimate \cite{cia}. 
 
For quasiregular mappings there is no 
way to patch together the argument shown in Subsection \ref{QG}. The key to making progress 
is to use the \em Stoilow decomposition \rm to translate the information we have  from the hypotheses 
into information about the analytic functions of the Stoilow decomposition. Let us recall the basics of the Stoilow decomposition, any $Q$-quasiregular 
mapping $u:\Omega\rightarrow \R^2$ can be written as the composition of a $Q$-quasiconformal homeomorphism $w_u:\Omega\rightarrow \R^2$ and an 
analytic function $\phi_u:w_u(\Omega)\rightarrow \R^2$ so that  
\begin{equation}
\label{introeq1}
u(z)=\phi_u(w_u(z)). 
\end{equation}
A good reference are the monographs of Astala-Iwaniec-Martin \cite{astala} Section 5.5.\ and Ahlfors \cite{alfors}.
 
The heart of the Stoilow decomposition is the fact that it is possible to solve Beltrami's equation. This allows us to find a $Q$-quasiconformal 
mapping $w_u$ that has the same \em Beltrami Coefficient  \rm as $Du$. The \em Beltrami Coefficient  \rm of a matrix $M$ is a $2\times 2$ conformal matrix $\mu_M$ (or more typically 
a complex number) that encodes the \em geometry \rm of the deformation of the unit ball by $M$, but \em not \rm the orientation 
or the size (formally $\lt[M\rt]_a \GI=\mu_{M} \lt[M\rt]_c$ where $\lt[M\rt]_c,\lt[M\rt]_a$ are the conformal and anti-conformal parts of $M$ and 
$\GI$ is a reflection across the $y$-axis, see Subsection \ref{conanticon} for more details). By solving Beltrami's equation we can find a homeomorphism $w_u$ with 
the property that 
\begin{equation}
\label{introeq2}
\mu_{D w_u(z)}=\mu_{D u(z)}\text{ for }a.e.\ z\in B_1 
\end{equation}
and  $w_u(z)-z=O(1/z)$, $w_v(z)-z=O(1/z)$. So for any $z\in B_1$ the shape of the image 
of the unit ball under $Du(z)$  is \em similar \rm to the shape of the image 
of the unit ball under $Dw_u(z)$. Hence the factorization represented by (\ref{introeq1}) is entirely natural.

Now the symmetric part of a gradient encodes both the geometry and the size. So a key result that starts 
the proof is a bound of the difference between Beltrami coefficients of two $Q$-quasiconformal matrices $A$, $B$ by $\lt|S(A)-S(B)\rt|$
\begin{equation}
\label{introeq3}
\lt|\mu_A-\mu_B\rt|\leq  32 \sqrt{Q}\min\lt\{\det(A)^{-\frac{1}{2}},\det(B)^{-\frac{1}{2}}\rt\}\lt|S(A)-S(B)\rt|.
\end{equation}
This is the contents of Lemma \ref{L.5}. Note as the determinants of $Q$-quasiconformal matrices $A, B$ 
get very small their norm gets very small so $\lt|S(A)-S(B)\rt|\leq \lt|S(A)\rt|+\lt|S(B)\rt|\leq 4 Q(\det(A)+\det(B))\approx 0$ but the \em geometry \rm of the deformation of the unit ball by $A$, $B$ could 
be very different hence the factor of $\min\lt\{\det(A)^{-\frac{1}{2}},\det(B)^{-\frac{1}{2}}\rt\}$ in the 
right hand side (\ref{introeq3}) is to be expected. 

Now the solutions of the Beltrami equation $w_u$, $w_v$ are essentially given by solving $\mathcal{C}(1-\mu_{Du} \mathcal{S})^{-1}$ and $\mathcal{C}(1-\mu_{Dv} \mathcal{S})^{-1}$ where $\mathcal{C}$ is the Cauchy transform and $\mathcal{S}$ is the Beurling-Ahlfors transform. Hence it should seem 
reasonable that we can prove an estimate showing $Dw_u$, $Dw_v$ are close in $L^p$ norm. As a consequence we establish 
\begin{equation}
\label{introeq6}
\|w_u-w_v\|_{L^{\infty}(B_{\frac{1}{2}})}\leq c \ep^{\frac{p}{120 Q^2}}. 
\end{equation}
This is part of the contents of Lemma \ref{LL1} and Lemma \ref{L.5}. 

Having established a quantitative relation between $w_u, w_v$ in order 
to prove the estimate on $Du$, $Dv$ we need to establish the relation $\phi_v'-\zeta\phi_u'\approx 0$ 
for some $\zeta\in \mathbb{C}$ with $\lt|\zeta\rt|=1$. We will establish this relation by applying Theorem \ref{MT} but first we have to set up some preliminary estimates. Since $w_u$ is a solution of the Beltrami equation 
we have explicit estimates on its $L^p$ norm and the $L^p$ norm of its inverse in terms 
of $Q$. Hence we are able to establish the existence of a constant $\mu=\mu(Q)$ such that 
\begin{equation}
\label{introeq11}
B_{2\mu}(w_u(0))\subset w_u(B_{\frac{1}{2}})\text{ and } B_{2\mu}(w_v(0))\subset w_v(B_{\frac{1}{2}}).
\end{equation}
This is the contents of part of Lemma \ref{LA3} and Lemma \ref{L1.4}. 

Now by (\ref{introeq6}) we know $w_u(B_{\frac{1}{2}})\subset 
w_v(B_1)$ so $\phi_u$ and $\phi_v$ are both defined on this set. Since the hypotheses are that the symmetric part of gradient are close we also know the size of the gradients $\na u$ and $\na v$ are close. 
By the chain rule this implies an estimate of the form 
\begin{equation}
\label{introeq11}
\int_{w_u(B_{\frac{1}{2}})} \lt|\lt|\phi_u'\rt|^2-\lt|\phi_v'\rt|^2\rt| dx\leq c_1 \ep^{\frac{p}{120 Q^2}}, 
\end{equation}
this is the content of Lemma \ref{L2}. We would 
like to apply Theorem \ref{MT} so a natural thing to do would be to use Cauchy's Theorem 
to find an analytic function $\psi$ such that $\psi'=\frac{\phi'_v}{\phi'_u}$ then establish 
appropriate lower bounds on $\lt|\phi'_u\rt|$ on some ball $B_{h_0}(x_0)$ to 
conclude 
\begin{equation}
\label{introeq12}
\int_{B_{h_0}(x_0)} \lt|1-\lt|\psi'(z)\rt|\rt|^2 dz\leq c_2 \ep^{\frac{p}{120 Q^2}}.
\end{equation}
 The non-degeneracy condition 
$\int_{B_1} \det(Du(z))^{-p} dz\leq C_p$ allows to find such a ball 
centred somewhere in $B_{\frac{\mu}{2}}(w_u(0))$, this is a the 
contents of Lemma \ref{L3}. Specifically we find some $h_0=h_0(Q,C_p)>0$ and some $\varpi=\varpi(Q,C_p)>0$ such 
that for some $x_0\in B_{\frac{\mu}{2}}(w_u(0))$ 
\begin{equation}
\label{introeq7}
\inf\lt\{\lt|\phi_u'(y)\rt|:y\in B_{h_0}(x_0)\rt\}\geq \varpi. 
\end{equation}

Let $\ti{\psi}(x,y)=(\mathrm{Re}(\psi(x+iy)),\mathrm{Im}(\psi(x+iy)))$. Reformulating (\ref{introeq12}) in 
matrix notation gives $\int_{B_{h_0}(x_0)} \mathrm{dist}^2(D\ti{\psi},SO(2)) dz\leq c_2 \ep^{\frac{p}{120 Q^2}}$. 
So we can apply Theorem \ref{MT}, however for reasons we will explain later we will instead use a more 
restricted version of it given by Proposition \ref{PP1} proved in Appendix. So we can conclude there 
exists some rotation $R$ such that 
\begin{equation}
\label{introeq13}
\int_{B_{h_0}(x_0)} \lt|D\ti{\psi}-R\rt| dz\leq c_3 \ep^{\frac{p}{960 Q^2}}.
\end{equation}
Returning this into complex notation and unwrapping it using the definition of $\psi$ we have 
\begin{equation}
\label{introeq14}
\int_{B_{h_0}(x_0)} \lt|\phi_v'(z)-\zeta \phi_u'(z)\rt| dz\leq c_4 \ep^{\frac{p}{960 Q^2}}.
\end{equation}

We need to extend control on $\phi'_v-\zeta \phi'_u$ to include an explicit neighbourhood of 
$w_u(0)$. We are able to do this by the fact that we are dealing with an analytic function $\phi_v-\zeta \phi_u$ and so have Talyor's Theorem. Since we already know $B_{2\mu}(w_u(0))\overset{(\ref{introeq11})}{\subset} w_u(B_1)$ and $x_0\in B_{\frac{\mu}{2}}(w_u(0))$ so we can use Talyor's theorem 
to extend control to $B_{\mu}(x_0)$ which contains $B_{\frac{\mu}{2}}(w_u(0))$

So let $w(z)=\phi_v'(z)-\zeta_1 \phi'_u(z)$. By the local Talyor Theorem we have 
$w(z)=\sum_{k=0}^{m} \frac{w^{(k)}(x_0)}{k!}(z-x_0)^k+(z-x_0)^{m+1}w_m(z)$ where $w_m(z)=\frac{1}{2\pi i}\int_{\partial B_{\frac{3\mu}{2}}(x_0)} \frac{w(\zeta)}{(\zeta-x_0)^m(\zeta-z)} d\zeta$. 

By the Coarea formula we can find $q\in (\frac{h_0}{8},h_0)$ such that $\int_{\partial B_q(x_0)} \lt|w(z)\rt| dH^1 z\leq 8 c_4  \ep^{\frac{p}{960 Q^2}}$. So by Cauchy's integral formula 
$$
\lt|w^{(k)}(x_0)\rt|=\frac{k!}{2\pi}\int_{\partial B_q(x_0)} \lt|\frac{w(\zeta)}{(\zeta-x_0)^{k+1}}\rt| d\zeta\leq 
\frac{4 c_4 k!}{\pi} \frac{\ep^{\frac{p}{960Q^2}}}{q^{k+1}}.
$$
We can also use the upper bound $\|Du\|_{L^2(B_1)}\leq 1$ and the upperbounds on $w_u$, $w_v$ to get upper 
bounds on $\phi_u$ and $\phi_v$ on $B_{\mu}(w_u(0))$ (this is part of the contents 
of Lemma \ref{L1.4}) so can estimate the remainder term 
$\|w_m\|_{L^{\infty}(B_{\mu}(z_0))}\leq 64\pi \mu^{-2}\lt(\frac{3\mu}{2}\rt)^{1-m}$. Thus  
we have 
\begin{equation}
\label{introeq19}
\lt|w(z)\rt|\leq \sum_{k=0}^m c_5 \ep^{\frac{p}{960 Q^2}} \lt(\frac{\mu}{q}\rt)^k
+64\pi\lt(\frac{3}{2}\rt)^{1-m}\text{ for any }z\in B_{\mu}(w_u(0)).
\end{equation}

The key is to make the right choice of $m$. If we choose $m$ too large then 
$\sum_{k=0}^m c_5 \lt(\frac{\mu}{q}\rt)^k$ will dominate $\ep^{\frac{p}{960 Q^2}}$ and 
the upperbound will be weak. If $m$ is too small then $64\pi\lt(\frac{3}{2}\rt)^{1-m}$ will 
not be small enough. The answer to to find $m$ that roughly equalizes these two 
quantities. An \em essential \rm point is that finding this $m$ requires knowing what the constants 
$h_0$, $c_5$, $\mu$ are. To estimate these constants we need to know $ c_1, c_2, c_3, c_4$ and $\varpi$ in 
 (\ref{introeq11}), (\ref{introeq12}), (\ref{introeq13}), (\ref{introeq14}) and (\ref{introeq7}). For this reason much effort will 
be made to track all constants in the estimates in this paper, since the methods are not close to 
being sharp we do not attempt to consistently calculate the best possible constants, but we do 
make efforts to prevent the constants blowing up too much throughout the paper. The reason we 
need the simplified version of Theorem \ref{MT} that is given by Proposition \ref{PP1} is that we need to know explicitly the constant in 
this inequality. This requires us to rewrite the proof of an estimate from \cite{fmul} while 
tracking the constants. The fact we are able to do this with the methods of \cite{fmul} is one of the reasons 
that Theorem \ref{T1} was not in practical terms accessible before the ideas introduced in \cite{fmul}. 
So making these estimates (recalling the fact $x_0\in B_{\frac{\mu}{2}}(w_u(0)))$ we have 
\begin{equation}
\label{introeq21}
\|\phi_v'-\zeta\phi_u'\|_{L^{\infty}(B_{\frac{\mu}{2}}(w_u(0)))}\leq c_5 C_p  \ep^{\frac{p^3}{4\times 10^7 Q^5 \log(10 C_p Q)}}
\end{equation}
This is the contents of Lemma \ref{L4}. 
By using the estimates on the closeness of $Dw_u$ and $Dw_v$ in $L^p$ we can then conclude 
that for some constant $\gamma=\gamma(Q)$ that 
\begin{equation}
\label{introeq21}
\|Dv-R Du\|_{L^{2}(B_{\gamma})}\leq c C_p \ep^{\frac{p^3}{4\times 10^7 Q^5 \log(10 C_p Q) }}.
\end{equation}
This is the contents of Proposition \ref{P1} below. Theorem \ref{T1} follows by a straightforward covering argument that gives estimate (\ref{fest1}).

\begin{a5}
\label{P1}
Suppose $v,u\in W^{1,2}(B_1:\R^2)$ 
are a $Q$-quasiregular mappings with $\int_{B_1} \det(Du)^{-p} dz\leq C_p$ for 
some $p\in (0,1)$ and 
$\int_{B_1} \lt|Du\rt|^2 dz\leq 1$. If   
\begin{equation}
\label{peq501}
\int_{B_1} \lt|S(Du)-S(Dv)\rt|^2 dz=\ep
\end{equation}
then there exists $R\in SO(2)$ and constant $\gamma=\gamma(Q)>0$ such that 
\begin{equation}
\label{peg502}
\int_{B_{\gamma}} \lt|Dv-R Du\rt| dz\leq c C_p \ep^{\frac{p^3}{4\times 10^{7} Q^5\log(10 C_p Q)}}. 
\end{equation}
\end{a5}

\em Remark. \rm We can assume with loss of generality 
\begin{equation}
\label{assump1}
u(0)=0,
\end{equation}
since if not the quasiregular mapping defined by $\ui(x)=u(x)-u(0)$ has this property. 

\section{Preliminaries}

\subsection{Conformal, Anti-conformal decomposition of $2\times 2$ matrices}
\label{conanticon}

Given $A\in M^{2\times 2}$ we can decompose $A$ into conformal and anti-conformal parts as 
follows 
\begin{equation}
\label{ffeqf11}
\lt(\begin{matrix} a_{11} & a_{12} \\
a_{21} & a_{22} \end{matrix}\rt)=\frac{1}{2}\lt(\begin{matrix} a_{11}+a_{22} & -(a_{21}-a_{12}) \\
a_{21}-a_{12} & a_{11}+a_{22} \end{matrix}\rt)+\frac{1}{2}\lt(\begin{matrix} a_{11}-a_{22} & a_{21}+a_{12} \\
a_{21}+a_{12} & -(a_{11}-a_{22}) \end{matrix}\rt).
\end{equation}
So for arbitrary matrix $A$ let 
\begin{equation}
\label{nicd1}
\lt[A\rt]_c:=\frac{1}{2}\lt(\begin{matrix} a_{11}+a_{22} & -(a_{21}-a_{12}) \\
a_{21}-a_{12} & a_{11}+a_{22} \end{matrix}\rt)\text{ and }
\lt[A\rt]_a:=\frac{1}{2}\lt(\begin{matrix} a_{11}-a_{22} & a_{21}+a_{12} \\
a_{21}+a_{12} & -(a_{11}-a_{22}) \end{matrix}\rt).
\end{equation}

It will often be convenient to write this decomposition as $A=\alpha R_{\theta}+\beta N_{\psi}$ 
where 
$$
R_{\theta}:= \lt(\begin{matrix} \cos{\theta} & -\sin{\theta} \\
\sin{\theta} & \cos{\theta} \end{matrix}\rt)\text{ and  }
N_{\psi}:= \lt(\begin{matrix} \cos{\psi} & \sin{\psi} \\
\sin{\psi} & -\cos{\psi} \end{matrix}\rt)
$$
 
Let $\GI:=\lt(\begin{smallmatrix} 1 & 0\\ 0 & -1\end{smallmatrix}\rt)$. 
The \bf Beltrami Coefficient  \rm of a matrix $A$ that relates the conformal and anti-conformal parts of $A$ 
is the conformal matrix $\mu_A$ defined by 
\begin{equation}
\label{eq2}
\lt[A\rt]_a \GI=\mu_{A} \lt[A\rt]_c.
\end{equation}

\begin{equation}
\label{eqcw100}
\|A\|^2\leq Q\det{A}\Rightarrow\frac{\lt(\alpha+\beta\rt)^2}{\alpha^2-\beta^2}\leq Q\Rightarrow \frac{\beta}{\alpha}
\leq\frac{Q-1}{Q+1}.
\end{equation}

Now notice that 
\begin{equation}
\label{eq5.2}
\lt|\mu_{A^{-1}}\rt|=\lt|\mu_A\rt|.
\end{equation}

And 
\begin{equation}
\label{nikb2}
\lt|\lt[A\rt]_a\rt|=\frac{1}{2}\lt|\lt(\begin{smallmatrix} a_{11} & a_{12}\\ a_{21} & a_{22}\end{smallmatrix}\rt)+\lt(\begin{smallmatrix} -a_{22} & a_{21}\\ a_{12} & -a_{11}\end{smallmatrix}\rt)\rt|\leq \lt|A\rt|.
\end{equation}

As $\beta N_{\psi} \GI=\mu_{A} \alpha R_{\theta}$, so 
\begin{equation}
\label{eqb3}
\lt|\mu_A\rt|=\sqrt{2}\frac{\beta}{\alpha}
\end{equation}

\subsection{The Beltrami equation}

The Beltrami equation is a linear complex PDE the relates the conformal part of the gradient to the 
anti-conformal, we briefly describe the connection between the classical complex formulation and 
and the matrix formulation we will be using in this paper. 

Take function from the complex plane to itself, $f(x+iy)=u(x,y)+iv(x,y)$. 
As is standard, $\frac{\partial}{\partial \bar{z}}f(x,y)=\frac{1}{2}(\partial_x+i\partial_y) f$ and 
$\frac{\partial}{\partial z} f(x,y)=\frac{1}{2}(\partial_x-i\partial_y) f$. 

 If we take a $\Omega\subset \mathbb{C}$ and a function $f:\Omega\rightarrow \mathbb{C}$ then define the $\R^2$ valued function $\ti{f}(x,y)=(\mathrm{Re}(f(x+iy)),\mathrm{Im}(f(x+iy)))$. Let $CO_{+}(2)$ denote the 
set of conformal $2\times 2$ matrices. And let $\lt[\cdot\rt]_M$ denote the homomorphism between $\mathbb{C}$ and $CO_{+}(2)$, so 
$\lt[a+ib\rt]_M=\lt(\begin{matrix} a & -b \\ b & a\end{matrix}\rt)$. 

So note 
\begin{equation}
\label{olleq6}
\sqrt{2}\lt|a+ib\rt|=\lt|\lt[a+ib\rt]_M\rt|.
\end{equation}
It is straight forward to see that 
$$
\lt[\frac{\partial f}{\partial z}\rt]_M=\lt[D\ti{f}\rt]_c\text{ and  }\lt[\frac{\partial f}{\partial \bar{z}}\rt]_M \lt(\begin{matrix} 1 & 0 \\ 0 & -1\end{matrix}\rt)=\lt[D\ti{f}\rt]_a,
$$
(recall the decomposition into conformal and anticonformal parts given by (\ref{ffeqf11}), (\ref{nicd1})).

Now as in 2.9.1.\ \cite{astala} letting $Df(z):\mathbb{C}\rightarrow \mathbb{C}$ denote the linear map that is the 
derivative of $f$ at $z$, then we have 
$Df(z)h=\frac{\partial f}{\partial z}(z)h+\frac{\partial f}{\partial \bar{z}}(z)\bar{h}$. Let $\lt[\cdot\rt]_{\mathbb{C}}$ be the identification of $\R^2$ with $\mathbb{C}$, i.e.\ $\lt[\lt(\begin{smallmatrix} a \\ b\end{smallmatrix}\rt)\rt]_{\mathbb{C}}=a+ib$. Let $f=u+iv$ so we have 
\begin{eqnarray}
Df(z)h&=&\frac{1}{2}\lt((u_x+v_y)+i(v_x-u_y)\rt)(h_1+i h_2)+\frac{1}{2}\lt(u_x-v_y+i(v_x+u_y)\rt)(h_1-i h_2)\nn\\
&=&\lt[\frac{1}{2}\lt(\begin{matrix} u_x+v_y & -(v_x-u_y) \\ v_x-u_y & u_x+v_y\end{matrix}\rt) \lt(\begin{matrix} h_1 \\ h_2\end{matrix}\rt)\rt]_{\mathbb{C}}+\lt[\frac{1}{2}\lt(\begin{matrix} u_x-v_y & v_x+u_y \\ v_x+u_y & -(u_x-v_y)\end{matrix}\rt) \lt(\begin{matrix} h_1 \\ h_2\end{matrix}\rt)\rt]_{\mathbb{C}}\nn\\
 &=&\lt[\lt(\frac{1}{2}\lt[D \ti{f}(x,y)\rt]_c+\frac{1}{2}\lt[D \ti{f}(x,y)\rt]_a \rt) \lt(\begin{matrix} h_1 \\ h_2\end{matrix}\rt)   \rt]_{\mathbb{C}}.
\end{eqnarray}

Given $f:\Omega\rightarrow \mathbb{C}$ one of the basic equations of Quasiregular analysis is the Beltrami equation
\begin{equation}
\label{peqq2}
\frac{\partial f}{\partial \bar{z}}(z)=\mu(z)\frac{\partial f}{\partial z}(z).
\end{equation}
As above define $\ti{f}=(\mathrm{Re}(f),\mathrm{Im}(f))$ then $\ti{f}$ satisfies 
\begin{equation}
\label{peqq3}
\lt[D \ti{f}(x,y)\rt]_a\II=\lt[\mu(x+iy)\rt]_M \lt[D \ti{f}(x,y)\rt]_c.
\end{equation}
By uniqueness this implies that 
\begin{equation}
\label{peqq4}
\lt[\mu(x+iy)\rt]_M=\mu_{D\ti{f}(x,y)}.
\end{equation}

The basic theorem about the solvability of the Beltrami equation (sometimes known as the measurable Riemann mapping theorem) is the following  
\begin{thm}[Morrey-Bojarski]
Suppose that $0\leq k<1$ and that $\lt|\mu(z)\rt|\leq \kappa\chara_{B_r}(z)$, $z\in\mathbb{C}$. Then 
there is a unique $f\in W^{1,p}_{loc}(\mathbb{C})$ (for every $p\in\lt[2,1+\frac{1}{\kappa}\rt)$)such that 
\begin{eqnarray}
&~&\frac{\partial f}{\partial \overline{z}}=\mu(z)\frac{\partial f}{\partial z}\text{ for almost every }z\in\mathbb{C}\nn\\
&~&f(z)=z+O\lt(\frac{1}{z}\rt)\text{ as }z\rightarrow \infty.
\end{eqnarray}
\end{thm}


\begin{deff}
Given a $Q$-quasiregular mapping $u$ we say the pair $w_u:B_1\rightarrow \R^2$, $\phi:w_u(B_1)\rightarrow u(B_1)$ are 
the Stoilow decomposition of $u$ iff 
\begin{equation}
\label{ffeqf5}
u(z)=\phi_u(w_u(z))\text{ for all }z\in B_1.
\end{equation}
Function $w_u$ is a $Q$-quasiregular mapping obtained by solving the Beltrami equation 
\begin{equation}
\label{peq70}
\lt[Dw_u(z)\rt]_{a}\II=\mu(z)\lt[Dw_u(z)\rt]_{c}
\end{equation}
where 
\begin{equation}
\label{peq71}
\mu(z)=\lt\{\begin{array}{lcl} \lt[Du(z)\rt]_a \II \lt[Du(z)\rt]_c^{-1} &\text{ for }& z\in B_1  \\
0 &\text{ for }&z\not \in B_1\end{array}\rt.. 
\end{equation}
\end{deff}
Note that (\ref{peq70}), (\ref{peq71}) are just the reformulation of the standard Beltrami equation 
and Beltrami coefficient in matrix notation as explained in Subsection \ref{conanticon} (\ref{eq2}) and equations (\ref{peqq3}), (\ref{peqq4}) of this subsection.   

As explained in the introduction, a consequence of (\ref{ffeqf5}), (\ref{peq70}) we have that 
$D\phi_u\in CO_{+}(2)=\lt\{\lm R: \lm>0, R\in SO(2) \rt\}$. So considered as a complex 
valued function of a complex variable, function $\phi_u$ is holomorphic. We will often consider $\phi_u$ as a 
holomorphic function of a complex variable without relabelling it.

\subsection{The Beltrami Coefficient  of gradient whose symmetric parts agree}

We require a Lemma 1 from \cite{lor27}. It is stated below  
\begin{a1} 
\label{L0.7}
Let $A\in M^{2\times 2}$, $\det(A)>0$. Let the Beltrami coefficient of $A$ be defined by (\ref{eq2}). The Beltrami coefficient  of $A$ and $A^{-1}$ are related in the following way 
\begin{equation}
\label{eqw1.1}
\mu_{A} \lt[A\rt]_c \GI=-\mu_{A^{-1}} \GI \lt[A\rt]_c.
\end{equation}
\end{a1}

\section{Lemmas for Theorem \ref{T1}}

\begin{a1} 
\label{L.5}
Suppose $A,B\in M^{2\times 2}$ with $\mathrm{det}(A)>0$ and $\mathrm{det}(B)>0$ and 
$\|A\|^2\leq Q\det(A)$, $\|B\|^2\leq Q\det(B)$ 
then
\begin{equation}
\label{ew1}
\lt|\mu_A-\mu_B\rt|\leq  \frac{32\sqrt{Q}}{\max\lt\{\sqrt{\det(A)},\sqrt{\det(B)}\rt\}}\lt|S(A)-S(B)\rt|.
\end{equation}
\end{a1}
\em Proof of Lemma \ref{L.5}. \rm 
Let $\|\|$ denote the operator norm. Since $\lt|A\rt|\leq \lt|A e_1\rt|+\lt|A e_2\rt|\leq 
2\|A\|$ we have the following estimate 
\begin{equation}
\label{nikb1}
\frac{\lt|A\rt|}{2}\leq \|A\|\leq \lt|A\rt|\text{ for any matrix }A\in M^{2\times 2}.
\end{equation}

Note $R(A)S(A)S(B)^{-1}R(B)^{-1}=AB^{-1}=\lt[AB^{-1}\rt]_a+\lt[AB^{-1}\rt]_c$. Thus 
\begin{equation}
\label{bbreq707}
S(A) S(B)^{-1}=R(A)^{-1}\lt[A B^{-1}\rt]_a R(B)+R(A)^{-1}\lt[A B^{-1}\rt]_c R(B)
\end{equation}
as the decomposition into conformal and anti-conformal parts are unique, so 
\begin{equation}
\label{bbreq704}
\lt[S(A)S(B)^{-1}\rt]_a\overset{(\ref{bbreq707})}{=}R(A)^{-1}\lt[A B^{-1}\rt]_a R(B).
\end{equation}
Note $\|ADJ(B)\|\overset{(\ref{nikb1})}{\leq} \lt|ADJ(B)\rt|=\lt|B\rt|\leq 2 \|B\|\leq 2 \sqrt{Q}\sqrt{\det(B)}$. 
So let $\delta=\lt|S(A)-S(B)\rt|$
\begin{eqnarray}
\label{nikd3}
\|S(A)S(B)^{-1}-Id\|&\leq& \|S(A)-S(B)\|\frac{\|ADJ(B)\|}{\det(B)}\nn\\
&\leq&\frac{2\|S(A)-S(B)\|\sqrt{Q}}{\sqrt{\det(B)}}\nn\\
&\overset{(\ref{nikb1})}{\leq}& \frac{2\delta \sqrt{Q}}{\sqrt{\det(B)}}. 
\end{eqnarray}
Now 
\begin{eqnarray}
\label{band03}
\lt|\lt[AB^{-1}\rt]_a\rt|&\overset{(\ref{bbreq704}),(\ref{nikb1})}{\leq}& 2\|\lt[S(A)S(B)^{-1}\rt]_a\|\nn\\
&=& 2\|\lt[S(A)S(B)^{-1}-Id\rt]_a\|\nn\\
&\overset{(\ref{nikb1})}{\leq}&2\lt|\lt[S(A)S(B)^{-1}-Id\rt]_a\rt|\nn\\
&\overset{(\ref{nikb2})}{\leq}&2\lt|S(A)S(B)^{-1}-Id\rt|\nn\\
&\overset{(\ref{nikb1}),(\ref{nikd3})}{\leq}&\frac{4\delta \sqrt{Q}}{\sqrt{\det(B)}}.
\end{eqnarray}
Thus as we know from (\ref{eqw1.1}) Lemma \ref{L0.7} applied to $B^{-1}$ that 
\begin{equation}
\label{peq97}
\mu_{B^{-1}}\lt[B^{-1}\rt]_c \II=-\mu_{B} \II \lt[B^{-1}\rt]_c,
\end{equation}
so 
\begin{eqnarray}
\label{bbreq703}
\lt[A B^{-1}\rt]_a&:=&\lt[\lt(\lt[A\rt]_c+\lt[A\rt]_a\rt)\lt(\lt[B^{-1}\rt]_c+\lt[B^{-1}\rt]_a\rt)\rt]_a\nn\\
&=&\lt[A\rt]_a \lt[B^{-1}\rt]_c+\lt[A\rt]_c \lt[B^{-1}\rt]_a\nn\\
&\overset{(\ref{eq2})}{=}&\mu_A \lt[A\rt]_c \II \lt[B^{-1}\rt]_c+\lt[A\rt]_c \mu_{B^{-1}} \lt[B^{-1}\rt]_c\II\nn\\
&\overset{(\ref{peq97})}{=}&\mu_A \lt[A\rt]_c \II \lt[B^{-1}\rt]_c-\lt[A\rt]_c \mu_B \II \lt[B^{-1}\rt]_c\nn\\
&=&(\mu_A-\mu_B)\lt[A\rt]_c \II \lt[B^{-1}\rt]_c.
\end{eqnarray}
For any matrix $A$ let $\Pi(A):=\inf\lt\{\lt|Av\rt|:\lt|v\rt|=1\rt\}$. Note that $\Pi(AB)\geq \Pi(A)\Pi(B)$.
Thus 
\begin{eqnarray}
\label{nikd5}
\Pi\lt(\mu_A-\mu_B\rt)\Pi\lt(\lt[A\rt]_c\rt)\Pi\lt(\lt[B^{-1}\rt]_c\rt)&\overset{(\ref{bbreq703})}{\leq}&\Pi\lt(\lt[AB^{-1}\rt]_a\rt)\leq \lt| \lt[AB^{-1}\rt]_a\rt|
\overset{(\ref{band03})}{\leq} \frac{4\delta \sqrt{Q}}{\sqrt{\det(B)}}.
\end{eqnarray}
Now $\Pi\lt(\lt[A\rt]_c\rt)=\sqrt{\lt(\frac{\det\lt(\lt[A\rt]_c\rt)}{\pi}\rt)}\geq \frac{\sqrt{\det(A)}}{2}$. And $\Pi\lt( \lt[B^{-1}\rt]_c\rt)\geq \frac{\sqrt{\det(B^{-1})}}{2}\geq \frac{1}{2\sqrt{\det(B)}}$. So putting these things together we have that 
$$
\frac{1}{4}\frac{\sqrt{\det(A)}}{\sqrt{\det(B)}}\Pi(\mu_A-\mu_B)\overset{(\ref{nikd5})}{\leq} \frac{4\delta \sqrt{Q}}{\sqrt{\det(B)}}.
$$
So $\Pi\lt(\mu_A-\mu_B\rt)\leq \frac{16\delta \sqrt{Q}}{\sqrt{\det(A)}}$.  By definition of $\Pi$ for any $\ep>0$ we can find 
$w\in S^1$ such that $\lt|(\mu_A-\mu_B)w\rt|\leq \frac{16\delta \sqrt{Q}}{\sqrt{\det(A)}}+\ep$. Since $\mu_A-\mu_B$ is 
conformal so $\lt|(\mu_A-\mu_B)e_1\rt|\leq \frac{16\delta \sqrt{Q}}{\sqrt{\det(A)}}+\ep$ and $\lt|(\mu_A-\mu_B)e_2\rt|\leq \frac{16\delta \sqrt{Q}}{\sqrt{\det(A)}}+\ep$ and thus  $\lt|\mu_A-\mu_B\rt|\leq \frac{32\delta \sqrt{Q}}{\sqrt{\det(A)}}$. Now since the hypotheses on $A$, $B$ are the same this implies 
$\lt|\mu_A-\mu_B\rt|\leq \frac{32\delta \sqrt{Q}}{\sqrt{\det(B)}}$ and hence we have established (\ref{ew1}). $\Box$

%
%

\subsection{Estimates on Beltrami equations} 

\subsubsection{Estimates of the Holder norm of solutions of the Beltrami equation}

We need bounds on the Holder norm of solutions of the Beltrami equation.


The first is a well known lemma whose constant we explicitly estimate. 

\begin{a1}
\label{LA1} Suppose $p>2$ and $u\in W^{1,p}(B_{w}(\zeta))$, then for any $x,y$ with 
\begin{equation}
\label{nikb4}
\lt|x-y\rt|<\frac{1}{2}\min\lt\{d(x,\partial B_{w}(\zeta))),d(y,\partial B_{w}(\zeta))) \rt\}
\end{equation}
\begin{equation}
\label{exc1}
\lt|u(x)-u(y)\rt|\leq 8\lt(\frac{p-1}{p-2}\rt)\lt|x-y\rt|^{\frac{p-2}{p}}\lt(\int_{B_{2\lt|x-y\rt|}(x)} \lt|\na u\rt|^p dx\rt)^{\frac{1}{p}}.
\end{equation}
\end{a1}
\em Proof of Lemma \ref{LA1}. \rm We will use the following Poincare type inequality (see page 267 \cite{evans2})
\begin{equation}
\label{awwaa1}
\int_{B_r(x)} \lt|u(x)-u(z)\rt| dz\leq \frac{r^2}{2}\int_{B_{r}(x)} \frac{\lt|\na u(z)\rt|}{\lt|z-x\rt|} dz
\end{equation}
Let $W=B_r(x)\cap B_r(y)$ with $r=\lt|x-y\rt|$. Note by (\ref{nikb4}), 
$B_{2r}(x)\subset B_{w}(\zeta)$. Let $p'$ denote the 
Holder conjugate of $p$. So 
\begin{eqnarray}
\label{nikb30.6}
\lt|u(x)-u(y)\rt|&\leq&\Xint{-}_{W} \lt|u(x)-u(z)\rt| dz+\Xint{-}_{W} \lt|u(y)-u(z)\rt| dz\nn\\
&\leq& \lt(\pi \lt(\frac{r}{2}\rt)^2\rt)^{-1}\lt(\int_{B_r(x)} \lt|u(x)-u(z)\rt| dz+\int_{B_r(y)} \lt|u(y)-u(z)\rt| dz   \rt)\nn\\ 
&\overset{(\ref{awwaa1})}{\leq}&\frac{2}{\pi}\int_{B_r(x)} \frac{\lt|\na u(z)\rt|}{\lt|x-z\rt|} dz
+\frac{2}{\pi}\int_{B_r(y)} \frac{\lt|\na u(z)\rt|}{\lt|y-z\rt|} dz\nn\\
&\leq&\frac{2}{\pi}\lt(\int_{B_r(x)} \lt|\na u\rt|^p \rt)^{\frac{1}{p}}
\lt(\int_{B_r(x)} \lt|x-z\rt|^{-p'} dz \rt)^{\frac{1}{p'}}\nn\\
&~&\qd+
\frac{2}{\pi}\lt(\int_{B_r(y)} \lt|\na u\rt|^p \rt)^{\frac{1}{p}}
\lt(\int_{B_r(y)} \lt|y-z\rt|^{-p'} dz \rt)^{\frac{1}{p'}}.
\end{eqnarray} 
Now 
\begin{eqnarray}
\lt(\int_{B_r(y)} \lt|y-z\rt|^{-p'} \rt)^{\frac{1}{p'}}&=&\lt(\int_{0}^r 2\pi s^{1-p'} ds \rt)^{\frac{1}{p'}}\nn\\
&\leq& 2\pi \lt( \frac{r^{2-p'}}{2-p'}\rt)^{\frac{1}{p'}}\nn\\
&=&2\pi \lt(\frac{p-1}{p-2}\rt)^{\frac{p-1}{p}} r^{\frac{p-2}{p}}.
\end{eqnarray}
Putting this together with (\ref{nikb30.6}) we have  
\begin{eqnarray}
\label{nikb30}
\lt|u(x)-u(y)\rt|&\leq& 4\lt(\frac{p-1}{p-2}\rt)^{\frac{p-1}{p}} r^{\frac{p-2}{p}}\lt(\int_{B_r(y)} \lt|\na u\rt|^p \rt)^{\frac{1}{p}}+4\lt(\frac{p-1}{p-2}\rt)^{\frac{p-1}{p}} r^{\frac{p-2}{p}}\lt(\int_{B_r(x)} \lt|\na u\rt|^p \rt)^{\frac{1}{p}}\nn\\
&\leq& 8\lt(\frac{p-1}{p-2}\rt)^{\frac{p-1}{p}} r^{\frac{p-2}{p}}\lt(\int_{B_{2r}(x)} \lt|\na u\rt|^p \rt)^{\frac{1}{p}}\nn
\end{eqnarray}
and hence we have established (\ref{exc1}). $\Box$


\begin{a1}
\label{LA2}

Suppose $0\leq \kappa<1$ and $\mu:\R^2\rightarrow \mathbb{C}$ is measurable and for 
some $x_0\in \R^2$, $\lt|\mu(z)\rt|\leq \kappa \cha_{B_{\tau}(x_0)}(z)$ for all $z\in \R^2$ and $f$ is a principle 
solution of the Beltrami equation 
$$
\frac{\partial f}{\partial \overline{z}}(z)=\mu(z)\frac{\partial f}{\partial z}(z).
$$
Let $p\in (2,2+\frac{1-\kappa}{3\kappa})$. For any $x\in \R^2$, $r>0$ we have 
\begin{equation}
\label{eqc7}
\lt(\int_{B_r(x)} \lt|\na f\rt|^p dz\rt)^{\frac{1}{p}}\leq r^{\frac{2}{p}}+\frac{2(1+3(p-2))}
{1-\kappa(1+2(p-2))}\tau^{\frac{2}{p}}.
\end{equation}
\end{a1}
\em Proof of Lemma \ref{LA2}. \rm Let $\SI$ denote the Beurling transform, let $S_p$ denote the 
$L_p$ norm of $\SI$. Consider the operator 
$$
(Id-\mu \SI)^{-1}=Id+\mu\SI+\mu\SI\mu\SI+\mu\SI\mu\SI\mu\SI\dots.
$$
Note that if $\phi\in L^p$ then 
\begin{equation}
\label{breq975}
\|\mu\SI\mu\SI\dots \mu \SI \phi\|_{L^p(\mathbb{C})}\leq (\kappa S_p)^n \|\phi\|_{L^p(\mathbb{C})}. 
\end{equation}

So we require $\kappa S_p<1$ in order for $(Id-\mu \SI)^{-1}$ to be well defined. By inequality 
(4.89) Section 4.5.2 \cite{astala} we have 
\begin{equation}
\label{awa60}
S_p<1+3(p-2).
\end{equation}

Thus it is sufficient for $\kappa(1+3(p-2))<1$ which is equivalent to $p<\frac{1-\kappa}{3\kappa}+2$. If this inequality is 
satisfied then 
\begin{equation}
\label{bbreq1}
\|(Id-\mu \SI)^{-1}\phi\|_{L^p}\overset{(\ref{breq975}),(\ref{awa60})}{\leq} \sum_{m=0}^{\infty} \lt(\kappa(1+3(p-2))\rt)^m\|\phi\|_{L^p(\mathbb{C})}\leq \frac{1}{1-\kappa(1+3(p-2))}\|\phi\|_{L^p}.
\end{equation}
So defining $\sigma=\CI\lt(\lt(Id-\mu \SI\rt)^{-1}\mu\rt)$ where $\CI$ is the Cauchy transform. As in 
the proof of Theorem 5.1.1. \cite{astala} we know that 
\begin{equation}
\label{nike1}
\frac{\partial \sigma}{\partial \overline{z}}=\frac{\partial}{\partial \overline{z}}
\lt(\CI\lt(\lt(Id-\mu \SI\rt)^{-1}\mu\rt)\rt)=\lt(Id-\mu \SI\rt)^{-1}\mu
\end{equation}
and 
\begin{equation}
\label{nik1}
\frac{\partial \sigma}{\partial z}=\frac{\partial}{\partial z}
\lt(\CI\lt(\lt(Id-\mu \SI\rt)^{-1}\mu\rt)\rt)=\SI\lt(\lt(Id-\mu \SI\rt)^{-1}\mu\rt).
\end{equation}
Thus 
\begin{equation}
\label{brreq2}
\|\frac{\partial \sigma}{\partial \overline{z}}\|_{L^p(\mathbb{C})}\overset{(\ref{bbreq1}),(\ref{nike1})}{\leq} 
\frac{\tau^{\frac{2}{p}}}{1-\kappa(1+3(p-2))}
\end{equation}
and 
\begin{equation}
\label{nike2}
\|\frac{\partial \sigma}{\partial z}\|_{L^p(\mathbb{C})}\overset{(\ref{nik1}),(\ref{bbreq1})}{\leq} 
S_p\frac{\tau^{\frac{2}{p}}}{1-\kappa(1+3(p-2))}\overset{(\ref{awa60})}{\leq} \frac{(1+3(p-2))\tau^{\frac{2}{p}}}{1-\kappa(1+3(p-2))}.
\end{equation}
Hence 
\begin{equation}
\label{eqc2}
\|\na \sigma\|_{L^p(\mathbb{C})}\overset{(\ref{brreq2}),(\ref{nike2})}{\leq} \frac{2(1+3(p-2))\tau^{\frac{2}{p}}}{1-\kappa(1+3(p-2))}.
\end{equation}

Now as in the proof of Theorem 5.1.2 \cite{astala} we see that $f(z)=z+\sigma(z)$ is 
the principle solution of the Beltrami 
equation, i.e. the function that satisfies 
$$
\frac{\partial f}{\partial \overline{z}}(z)=\mu(z) \frac{\partial f}{\partial z}(z)\text{ for 
a.e.\ }z
$$
and $f(z)=z+\OI(\frac{1}{z})$ as $z\rightarrow \infty$. So note that for any $x\in \R^2$ we have that 
$$
\lt(\int_{B_r(x)} \lt|\na f\rt|^p dz\rt)^{\frac{1}{p}}\leq r^{\frac{2}{p}}
+\lt(\int_{B_r(x)} \lt|\na \sigma\rt|^p dz\rt)^{\frac{1}{p}}.
$$
Putting this together with (\ref{eqc2}) we have (\ref{eqc7}). $\Box$


\begin{a1}
\label{LA3}
Suppose $0 \leq \kappa<1$ and $\mu:\R^2\rightarrow \mathbb{C}$ is measurable and 
$\lt|\mu(z)\rt|\leq \kappa \cha_{B_1}(z)$ for all $z\in \R^2$. 
Let $f:\R^2\rightarrow \mathbb{C}$ be the principle solution of the Beltrami equation 
$$
\frac{\partial f}{\partial \overline{z}}(z)=\mu(z) \frac{\partial f}{\partial z}(z)\text{ for 
a.e.\ }z
$$
and let $h:\mathbb{C}\rightarrow \mathbb{C}$ be the global inverse of $f$. Let 
$p\in (2,2+\frac{1-\kappa}{3\kappa})$. Then 
\begin{equation}
\label{eqnzx1}
\lt|f(x_1)-f(x_2)\rt|\leq 48\lt|x_1-x_2\rt|^{\frac{p-2}{p}}
\lt(\frac{p-1}{p-2}\rt)\lt(\frac{1+3(p-2)}{1-\kappa(1+2(p-2))}\rt)\text{ for any }x_1,x_2\in B_1.
\end{equation}
And
\begin{equation}
\label{eqnzx2}
\lt|h(y_1)-h(y_2)\rt|\leq 2400\lt|y_1-y_2\rt|^{\frac{p-2}{p}}
\lt(\frac{p-1}{p-2}\rt)^2\lt(\frac{1+3(p-2)}{1-\kappa(1+2(p-2))}\rt)^2\text{ for any }y_1,y_2\in f(B_1).
\end{equation}
As a consequence for any $B_r(x)\subset B_1$ 
\begin{equation}
\label{eqe4.5}
B_{\lt(\frac{(1-\kappa)^{6} r}{4.3\times 10^{10}}\rt)^{\frac{12}{1-\kappa}}}(f(x))\subset f\lt(B_{r}(x)\rt).
\end{equation}
In addition for any $\alpha>0$ such that 
$B_{\lt(\frac{\alpha\lt(1-\kappa\rt)^{2}}{2304}\rt)^{\frac{12\kappa}{1-\kappa}}}(x)\subset B_1$
\begin{equation}
\label{ezq3.5}
f\lt(B_{\lt(\frac{\alpha\lt(1-\kappa\rt)^{2}}{3456}\rt)^{\frac{12}{1-\kappa}}}(x)\rt)\subset B_{\alpha}(f(x)). 
\end{equation}
\end{a1}

\em Proof of Lemma \ref{LA3}. \rm By Lemma \ref{LA2} we have that 
\begin{equation}
\label{eqnzx3}
\lt(\int_{B_4(x)} \lt|\na f\rt|^p dz \rt)^{\frac{1}{p}}\overset{(\ref{eqc7})}{\leq} 4+\frac{2(1+3(p-2))}{1-\kappa(1+2(p-2))}\text{ for any }x\in \R^2. 
\end{equation}
So by Lemma \ref{LA1} we know that 
\begin{eqnarray}
\label{ava350}
\lt|f(x_1)-f(x_2)\rt|&\overset{(\ref{exc1})}{\leq}& 8\lt(\frac{p-1}{p-2}\rt)\lt|x_1-x_2\rt|^{\frac{p-2}{p}}
\lt(\int_{B_{2\lt|x_1-x_2\rt|}(x)} \lt|\na f\rt|^p dz \rt)^{\frac{1}{p}}\nn\\
&\overset{(\ref{eqnzx3})}{\leq}&8\lt|x_1-x_2\rt|^{\frac{p-2}{p}}\lt(\frac{p-1}{p-2}\rt)
\lt(4+\frac{2(1+3(p-2))}{1-\kappa(1+2(p-2))} \rt)\nn\\
&\leq& 48\lt|x_1-x_2\rt|^{\frac{p-2}{p}}\lt(\frac{p-1}{p-2}\rt)
\lt(\frac{(1+3(p-2))}{1-\kappa(1+2(p-2))} \rt)
\end{eqnarray}
(using the fact $\kappa(1+2(p-2))\in (0,1)$ for the last inequality) so estimate (\ref{eqnzx1}) holds true. 

Now if we consider the Beltrami equation of $f$ we have 
$(\na f(x))_a \GI=\mu_{\na f(x)}(\na f(x))_c$, so if $z=f(x)$ then 
$\na h(z)=\lt(\na f(h(z))\rt)^{-1}=(Df(x))^{-1}$. Now the Beltrami equation for $h$ 
is $(\na h(z))_a \GI=\mu_{\na h(z)}(\na h(z))_c$. By (\ref{eq5.2}) we have that 
\begin{equation}
\label{peqq1.7}
\lt|\mu_{\na h(z)}\rt|=\lt|\mu_{\lt(\na f(x)\rt)^{-1}}\rt|=\lt|\mu_{\na f(x)}\rt|. 
\end{equation}
Now if $z\not \in f(B_1)$, since $\na h(z)=\lt(\na f(h(z))\rt)^{-1}$ and since $h(z)\not \in B_1$, $\na f(h(z))\in CO_{+}(2)$ so $\na h(z)\in CO_{+}(2)$ thus $\mu_{\na h(z)}=0$. 

Let 
\begin{equation}
\label{nikb40}
\Lambda_p^{\kappa}=48\lt(\frac{p-1}{p-2}\rt)\lt(\frac{1+3(p-2)}{1-\kappa(1+2(p-2))}\rt).
\end{equation}
Note that for any $x\in B_1$, $\lt|f(x)-f(0)\rt|\overset{(\ref{eqnzx1})}{\leq} \Lambda_p^{\kappa}$ so 
\begin{equation}
\label{nikb70}
f(B_1)\subset B_{ \Lambda_p^{\kappa}}(f(0)).
\end{equation}
So returning to complex notation we have $\frac{\partial h}{\partial \overline{w}}(w)=\gamma(w)
\frac{\partial h}{\partial w}(w)$ where 
$$
\lt|\gamma(w)\rt|\overset{(\ref{peqq4}),(\ref{peqq1.7})}{\leq} \kappa \cha_{B_{\Lambda_p^{\kappa}}(f(0))}.
$$ 
 
By Lemma \ref{LA2} (\ref{eqc7}) we know 
\begin{eqnarray}
\label{peq110}
\lt(\int_{B_{4 \Lambda_p^{\kappa}}(y_1)} 
\lt|\na h\rt|^p dz\rt)^{\frac{1}{p}}&\leq& \lt(4 \Lambda_p^{\kappa}\rt)^{\frac{2}{p}}+\lt(\frac{2(1+3(p-2))}{1-\kappa(1+2(p-2))}\rt)\lt(\Lambda_p^k\rt)^{\frac{2}{p}}\nn\\
&\overset{(\ref{nikb40})}{\leq}& 300\lt(\frac{p-1}{p-2}\rt)\lt(\frac{1+3(p-2)}{1-\kappa(1+2(p-2))}\rt)^2.
\end{eqnarray}
Now for any $y_1,y_2\in f(B_1)$ by (\ref{nikb70}) we know 
$y_2\in B_{2\Lambda_p^{\kappa}}(y_1)$ so by Lemma \ref{LA1}
\begin{equation}
\label{nikb50}
\lt|h(y_1)-h(y_2)\rt|\overset{(\ref{peq110}),(\ref{exc1})}{\leq} 2400 \lt(\frac{p-1}{p-2}\rt)^2
\lt|y_1-y_2\rt|^{\frac{p-2}{p}}\lt(\frac{1+3(p-2)}{1-\kappa(1+2(p-2))}\rt)^2
\end{equation}
and hence (\ref{eqnzx2}) is established.

Now suppose $B_r(x)\subset B_1$. Let 
\begin{equation}
\label{nikb7}
p=\min\lt\{2+\frac{1-\kappa}{6\kappa},3\rt\}. 
\end{equation}
If $p<3$ then 
\begin{eqnarray}
\label{eqe3}
 \lt(\frac{p-1}{p-2}\rt)^2\lt(\frac{1+3(p-2)}{1-\kappa(1+2(p-2))}\rt)^2&\overset{(\ref{nikb7})}{\leq}&
\lt(\frac{12}{1-\kappa}\rt)^2\lt(\frac{4}{1-\kappa\lt(1+2\lt(\frac{1-\kappa}{6\kappa}\rt)\rt)}\rt)^2\nn\\
&=& \lt(\frac{12}{1-\kappa}\rt)^2 \lt(\frac{6}{1-\kappa}\rt)^2\nn\\
&=& \frac{(72)^2}{(1-\kappa)^4}.
\end{eqnarray}
If $p=3$ then $2+\frac{1-\kappa}{6\kappa}\geq 3$, $1-\kappa\geq 6\kappa$ so $0<\kappa\leq \frac{1}{7}$. So 
\begin{eqnarray}
\label{eqe4}
 \lt(\frac{p-1}{p-2}\rt)^2\lt(\frac{1+3(p-2)}{1-\kappa(1+2(p-2))}\rt)^2&=& 
4\lt(\frac{4}{1-3\kappa}\rt)^2\nn\\
&\leq&  4\lt( \frac{4}{1-\frac{3}{7}}\rt)^2\nn\\
&\leq& \frac{200}{(1-\kappa)^4}.\nn
\end{eqnarray}
So for any $p$ we have that 
\begin{equation}
\label{bbreq3}
\lt(\frac{p-1}{p-2}\rt)^2\lt(\frac{1+3(p-2)}{1-\kappa(1+2(p-2))}\rt)^2\leq \frac{(72)^2}{(1-\kappa)^4}
\end{equation}
and so 
\begin{equation}
\label{nikb9}
\frac{\Lambda_p^{\kappa}}{48}\overset{(\ref{nikb40})}{=}\lt(\frac{p-1}{p-2}\rt)\lt(\frac{1+3(p-2)}{1-\kappa(1+2(p-2))}\rt)\leq \frac{72}{(1-\kappa)^{2}}.
\end{equation}
Thus 
\begin{equation}
\label{nikc1.5}
\Lambda_p^{\kappa}\leq \frac{3456}{(1-\kappa)^{2}}.
\end{equation}
Now if $p\in (2,3)$ then $2+\frac{1-\kappa}{6\kappa}<3$ so $\frac{1}{7}<\kappa<1$
\begin{equation}
\label{nikd1}
\frac{p-2}{p}\overset{(\ref{nikb7})}{=}\frac{1-\kappa}{6p\kappa}
\overset{(\ref{nikb7})}{=}\frac{1-\kappa}{6\lt(2+\frac{1-\kappa}{6\kappa}\rt)\kappa}=
\frac{1-\kappa}{11\kappa+1}\geq\frac{1-\kappa}{12}
\end{equation}
and 
\begin{eqnarray}
\label{peq3}
\frac{p-2}{p}-\frac{1-\kappa}{12}&\overset{(\ref{nikd1})}{=}&\frac{1-\kappa}{11\kappa+1}-\frac{1-\kappa}{12}\nn\\
&<&\frac{1-\kappa}{12\kappa}-\frac{1-\kappa}{12}=
\frac{1-2\kappa+\kappa^2}{12\kappa}\nn\\
&<&\frac{1}{12\kappa}<\frac{7}{12}.
\end{eqnarray}
And if $p=3$ since $\kappa\leq \frac{1}{7}$ we have 
\begin{equation}
\label{peq1}
\frac{p-2}{p}=\frac{1}{3}\geq \frac{1-\kappa}{12}
\end{equation}
and 
\begin{equation}
\label{peq2}
0\leq \frac{p-2}{p}-\frac{1-\kappa}{12}\leq \frac{1}{3}<1.
\end{equation}
Thus in all cases we have 
\begin{equation}
\label{bbreq5}
\frac{p-2}{p}\overset{(\ref{peq1}),(\ref{nikd1})}{\geq} \frac{1-\kappa}{12}
\end{equation}
and 
\begin{equation}
\label{nikzz2}
0\leq \frac{p-2}{p}-\frac{1-\kappa}{12}\overset{(\ref{peq3}),(\ref{peq2})}{<}1.
\end{equation}
Let $\varpi\in f(\partial B_{r}(x))$ be such that $\lt|\varpi-f(x)\rt|=\inf\lt\{\lt|f(z)-f(x)\rt|:z\in \partial B_{r}(x)\rt\}$. So 
\begin{equation}
\label{nikc1}
\lt|\varpi-f(x)\rt|\overset{(\ref{eqnzx1})}{\leq} 48 r^{\frac{p-2}{p}}\lt(\frac{p-1}{p-2}\rt)\lt(\frac{1+3(p-2)}{1-\kappa\lt(1+2\lt(\frac{1-\kappa}{6\kappa}\rt)\rt)}\rt)\overset{(\ref{nikb40})}{=}\Lambda^{\kappa}_p r^{\frac{p-2}{p}}.
\end{equation}
So as 
\begin{eqnarray}
\label{peq4}
\lt|\varpi-f(x)\rt|^{\frac{p-2}{2}}&=&\lt|\varpi-f(x)\rt|^{\frac{1-\kappa}{12\kappa}}
\lt|\varpi-f(x)\rt|^{\frac{p-2}{p}-\frac{1-\kappa}{12\kappa}}\nn\\
&\overset{(\ref{nikzz2}),(\ref{nikc1})}{\leq}&\lt|\varpi-f(x)\rt|^{\frac{1-\kappa}{12\kappa}} \Lambda^{\kappa}_p,
\end{eqnarray}
thus 
\begin{eqnarray}
r&=&\lt|h(\varpi)-x\rt|\nn\\
&=&\lt|h(\varpi)-h(f(x))\rt|\nn\\
&\overset{(\ref{nikb50}),(\ref{bbreq3})}{\leq}& 2400\lt|\varpi-f(x)\rt|^{\frac{p-2}{p}}\times \frac{(72)^2}{(1-\kappa)^4}\nn\\
&\overset{(\ref{peq4})}{\leq}&\frac{2400\times (72)^2 \Lambda_p^{\kappa}}{(1-\kappa)^4}
\lt|\varpi-f(x)\rt|^{\frac{1-\kappa}{12}}\nn\\
&\overset{(\ref{nikc1.5})}{\leq}&\frac{4.3\times 10^{10}}{(1-\kappa)^{6}}\lt|\varpi-f(x)\rt|^{\frac{1-\kappa}{12}}.\nn
\end{eqnarray} 
Thus $\frac{(1-\kappa)^{6} r}{4.3\times 10^{10}}\leq \lt|\varpi-f(x)\rt|^{\frac{1-\kappa}{12}}$ so 
$\lt(\frac{(1-\kappa)^{6} r}{4.3\times 10^{10}}\rt)^{\frac{12}{1-\kappa}}\leq \lt|\varpi-f(x)\rt|$ which implies (\ref{eqe4.5}). Now finally 
\begin{eqnarray}
\lt|f(x)-f(y)\rt|&\overset{(\ref{eqnzx1}),(\ref{nikb40})}{\leq}&\Lambda^{\kappa}_p \lt|x-y\rt|^{\frac{p-2}{p}}\nn\\
&\overset{(\ref{nikc1.5})}{\leq}&\frac{3456}{(1-\kappa)^2} \lt|x-y\rt|^{\frac{p-2}{p}}\nn\\
&\overset{(\ref{bbreq5})}{\leq}&\frac{3456}{(1-\kappa)^2} \lt|x-y\rt|^{\frac{1-\kappa}{12}}\text{ for any }\lt|x-y\rt|<1.
\end{eqnarray}
Thus for any $y\in B_{\lt(\frac{\alpha(1-\kappa)^{2}}{3456}\rt)^{\frac{12}{1-\kappa}}}(x)$ we have 
$$
\lt|f(x)-f(y)\rt|< \alpha
$$ 
which implies (\ref{ezq3.5}). $\Box$

%
%
%
%

\begin{a1}
\label{L1.4}
For $\mu=\mucon$, $\gamma=\gammacon$ we have 
\begin{equation}
\label{eqz11}
B_{2\mu}(w_u(0))\subset w_u\lt(B_{\frac{1}{2}}\rt),\;\; B_{2\mu}(w_v(0))\subset w_v\lt(B_{\frac{1}{2}}\rt) 
\end{equation}
and 
\begin{equation}
\label{ezq1}
w_u(B_{\gamma})\subset B_{\frac{\mu}{2}}(w_u(0)),\;\; w_v(B_{\gamma})\subset B_{\frac{\mu}{2}}(w_v(0)).
\end{equation}
In addition  
\begin{equation}
\label{eq1}
\|\phi_u\|_{L^{\infty}(B_{2\mu}(w_u(0)))}\leq 4,\;\; \|\phi_v\|_{L^{\infty}(B_{2\mu}(w_v(0)))}\leq 4.
\end{equation}
\begin{equation}
\label{eq3}
\|\phi_u^{'}\|_{L^{\infty}(B_{\mu}(w_u(0)))}\leq \frac{16\pi}{\mu},\;\; \|\phi_v^{'}\|_{L^{\infty}(B_{\mu}(w_u(0)))}\leq \frac{16 \pi}{\mu}.
\end{equation}
and 
\begin{equation}
\label{eq5}
\|\phi_u^{''}\|_{L^{\infty}(B_{\mu}(w_u(0)))}\leq \frac{48\pi}{\mu^2},\;\; \|\phi_u^{''}\|_{L^{\infty}(B_{\mu}(w_u(0)))}\leq \frac{48\pi}{\mu^2}.
\end{equation}
\end{a1}

\em Proof of Lemma \ref{L1.4}. \rm We will argue the estimate for $u$, $\phi_u$. The estimates for $v$, $\phi_v$ follow by exactly the same arguments.

Now recall from (\ref{eqcw100}), (\ref{eqb3}) we can take 
\begin{equation}
\label{breq62}
\kappa=\frac{Q-1}{Q+1}. 
\end{equation}
Now 
\begin{equation}
\label{nikc3}
1-\kappa= \frac{Q+1-(Q-1)}{Q+1}= \frac{2}{Q+1}.
\end{equation}
Thus
\begin{eqnarray}
\lt(\frac{r}{10^{10} (Q+1)^6}\rt)^{6(Q+1)}&\leq& \lt( \frac{64}{10}\times\frac{r}{10^{10} (Q+1)^{6}}\rt)
^{6(Q+1)}\nn\\
&\overset{(\ref{nikc3})}{=}&
\lt(\frac{r 2^6}{10^{11} (Q+1)^{6}}\rt)^{\frac{12}{1-\kappa}}\nn\\
&\overset{(\ref{nikc3})}{=}&\lt(\frac{r (1-\kappa)^{6}}{10^{11}}\rt)^{\frac{12}{1-\kappa}}\nn\\
&\leq& \lt(\frac{r}{2}\frac{(1-\kappa)^{6}}{4.3\times 10^{10}}\rt)^{\frac{12}{1-\kappa}}\nn
\end{eqnarray}
So by Lemma \ref{LA3}
\begin{equation}
\label{nikf1}
B_{\lt(\frac{r}{10^{10}(Q+1)^6}\rt)^{6(Q+1)}}(w_u(x))\overset{(\ref{eqe4.5})}{\subset} w_u(B_{\frac{r}{2}}(x))
\end{equation}
so defining $\mu=\mucon$, since $2\mu\leq (10^{10}(Q+1)^6)^{-6(Q+1)}$ we do indeed have 
$B_{2\mu}(w_u(0))\overset{(\ref{nikf1})}{\subset} w_u(B_{\frac{1}{2}})$ and so 
(\ref{eqz11}) is established. 

Note we have 
$$
B_{\lt(\frac{\mu}{2000(Q+1)^2}\rt)^{6(Q+1)}}(0)\subset B_{\lt(\frac{4\mu}{2\times 3456(Q+1)^2}\rt)^{6(Q+1)}}(0)
\overset{(\ref{nikc3})}{=} 
B_{\lt(\frac{(1-\kappa)^{2}\mu}{2\times 3456}\rt)^{\frac{12}{1-\kappa}}}(0).
$$
So by (\ref{ezq3.5}) we have $w_u\lt(B_{\lt(\frac{\mu}{2000(Q+1)^2}\rt)^{6(Q+1)}}(0)\rt)\subset B_{\frac{\mu}{2}}(w_u(0))$ which establishes (\ref{ezq1}). 

As $\int_{B_1} \lt|Du\rt| dz\leq \sqrt{\pi}$ we can find $h\in (\frac{1}{2},1)$ such that 
$\int_{\partial B_h} \lt|Du\rt| dH^1 z\leq 2\sqrt{\pi}$. Since $u$ is open $\partial u(B_h)\subset u(\partial B_h)$ so 
$H^1(\partial u(B_h))\leq 4$. 

So
\begin{equation}
\label{peq101}
u(B_{\frac{1}{2}})\subset u(B_h)\subset B_4(u(0)).
\end{equation}
Now by (\ref{eqz11}) since $w_u$ is a homeomorphism $w_u^{-1}(B_{2\mu}(w_u(0)))\subset B_{\frac{1}{2}}$ so as $\phi_u=u\circ w_u^{-1}$ we have 
$$
\phi_u(B_{2\mu}(w_u(0)))\overset{(\ref{peq101})}{\subset} B_4(u(0))\overset{(\ref{assump1})}{=} B_4(0). 
$$
So
\begin{equation}
\label{nikc5}
\|\phi_u\|_{L^{\infty}(B_{2\mu}(w_u(0)))}\leq 4.
\end{equation}
Thus for any $z\in B_{\mu}(w_u(0))$
\begin{eqnarray}
\lt|\phi_u'(z)\rt|&\leq&\int_{\partial B_{2\mu}(w_u(0))} \frac{\lt|\phi_u(\zeta)\rt|}{\lt|\zeta-z\rt|^2} \lt|d\zeta\rt|\leq\sup_{B_{2\mu}(w_u(0))} \lt|\phi_u(\zeta)\rt|\int_{\partial B_{2\mu}(w_u(0))} \frac{1}{\lt|\zeta-z\rt|^2} \lt|d\zeta\rt|\leq \frac{16\pi}{\mu}.\nn 
\end{eqnarray}
In the same way for any $z\in B_{\mu}(w_u(0))$
\begin{eqnarray}
\lt|\phi_u^{''}(z)\rt|&\leq&3\sup_{B_{2\mu}(w_u(0))} \lt|\phi_u(\zeta)\rt|\int_{\partial B_{2\mu}(w_u(0))} \frac{1}{\lt|\zeta-z\rt|^3} \lt|d\zeta\rt|\leq \frac{48 \pi}{\mu^2}.\;\;\;\Box\nn 
\end{eqnarray}

%
%
%
%

\begin{a1}
\label{LL1}
Let $u\in W^{1,2}(B_1,\R^2)$, $v\in W^{1,2}(B_1,\R^2)$ be $Q$-regular functions. Suppose for some $p\in (0,1)$
\begin{equation}
\label{equj1} 
\int_{B_1} \det(Du)^{-p} dz\leq C_p
\end{equation}
and  
\begin{equation}
\label{equl1} 
\| S(\na u)-S(\na v)\|_{L^2(B_1)}=\ep^2.
\end{equation}

Let $w_u$, $w_v$ be the quasiconformal mappings we obtain from the Stoilow decomposition of $u$ and $v$ we have 
\begin{equation}
\label{equa127}
\|\na w_u-\na w_v\|_{L^2(B_1)}\leq 24\pi C_p \sqrt{Q}\ep^{\frac{p}{48 Q^2}}.
\end{equation}
\end{a1}

\em Proof of Lemma \ref{LL1}. \rm We will require Lemma 5.3.1 \cite{astala}. This lemma controls the $L^p$ difference 
between the solutions $f,g$ of the Beltrami equations  
$\frac{\partial f}{\partial \overline{z}}=\mu(z)\frac{\partial f}{\partial z}$, 
$\frac{\partial g}{\partial \overline{z}}=\nu(z)\frac{\partial g}{\partial z}$ where $\lt|\mu\rt|,\lt|\nu\rt|\leq \kappa\chara _{B_r}$. 
Specifically for $p\in \lt[2,1+\frac{1}{\kappa}\rt)$, (see p163 \cite{astala}) Lemma 5.3.1 asserts that 
$\|\frac{\partial f}{\partial \overline{z}}-  \frac{\partial f}{\partial \overline{z}}\|_{L^p(\mathbb{C})}\leq 
\|\mu-\nu\|_{L^{\frac{ps}{s-1}}(\mathbb{C})}$ where $s$ is a number such that $p<sp<1+\frac{1}{k}$. 

Recall from (\ref{eqcw100}), (\ref{eqb3}) we can take $\kappa=\frac{Q-1}{Q+1}$. So 
$\frac{1}{2}(3+\frac{1}{\kappa})=\frac{1}{2}(3+\frac{Q+1}{Q-1})<1+\frac{1}{\kappa}$. So define  
\begin{equation}
\label{eqvuu2}
\PPI_Q=\lt\{\begin{array}{ll} 2^{-1}(3+\frac{Q+1}{Q-1})  & \text{ for }Q\geq 2 \\
3 & \text{ for }Q< 2\end{array}\rt.
\end{equation}

Now from (\ref{eqcw100}) and (\ref{eqb3}) we have that $\lt|\mu_{\na w_u}\rt|<\sqrt{2}\lt(\frac{Q-1}{Q+1}\rt)$ and 
$\lt|\mu_{\na w_v}\rt|<\sqrt{2}\lt(\frac{Q-1}{Q+1}\rt)$.  We require $s>1$ to be such that 
$s\PPI_Q=\frac{s}{2}(3+\frac{Q+1}{Q-1})<1+\frac{Q+1}{Q-1}$, i.e.\ 
$s<\frac{2+2\frac{Q+1}{Q-1}}{3+\frac{Q+1}{Q-1}}$.

Define 
\begin{equation}
\label{eqvuu1}
s_Q=\lt\{\begin{array}{ll} \frac{2+2\frac{Q+1}{Q-1}}{\frac{5}{2}+\frac{3}{2}\frac{Q+1}{Q-1}} & \text{ for }Q\geq 2 \\
\frac{8}{7} & \text{ for }Q< 2.\end{array}\rt.
\end{equation}

Note if $Q\geq 2$, 
\begin{eqnarray}
\label{nikc11}
s \PPI_Q&=&\frac{1+\frac{Q+1}{Q-1}}{\frac{5}{2}+\frac{3}{2}\frac{Q+1}{Q-1}}\lt(\frac{3(Q-1)+Q+1}{Q-1}\rt)\nn\\
&=&\frac{\frac{2Q}{Q-1}}{\frac{5(Q-1)+3(Q+1)}{2(Q-1)}}\lt(\frac{4Q-2}{Q-1}\rt)=\lt(\frac{2Q}{4Q-1}\rt)\lt(\frac{4Q-2}{Q-1}\rt).
\end{eqnarray}
Let 
$$
\BI=\lt\{z\in B_1: \det(Du(z))\leq \sqrt{\ep}\rt\}.
$$
So 
\begin{equation}
C_p\geq \int_{B_1} \det(Du(z))^{-p} dz\geq \ep^{-\frac{p}{2}}\lt|\BI\rt|.\nn
\end{equation}
Thus 
\begin{equation}
\label{equv2}
\lt|\BI\rt|\leq C_p \ep^{\frac{p}{2}}.
\end{equation}
Now for any $z\in B_1\backslash \BI$ by Lemma \ref{L.5} we have 
\begin{eqnarray}
\label{equv23}
\lt|\mu_{Dw_u(z)}-\mu_{Dw_v(z)}\rt|&=&\lt|\mu_{Dv(z)}-\mu_{Du(z)}\rt|\overset{(\ref{ew1})}{\leq}\frac{32\sqrt{Q}}{\sqrt{\ep}}\lt|S(Du(z))-S(Dv(z))\rt|.
\end{eqnarray}
And note 
\begin{equation}
\label{equv3}
\lt|\mu_{Dv(z)}-\mu_{Du(z)}\rt|\leq 3\text{ for any }z\in B_1.
\end{equation}
Hence  
\begin{eqnarray}
\label{equv24}
\int_{B_1} \lt|\mu_{Dw_u(z)}-\mu_{Dw_v(z)}\rt| dz&\overset{(\ref{equv3}),(\ref{equv23})}{\leq}& 3\lt|\BI\rt|
+\frac{32\sqrt{Q}}{\sqrt{\ep}}\int_{B_1} \lt|S(Du(z))-S(Dv(z))\rt| dz\nn\\
&\overset{(\ref{equl1}),(\ref{equv2})}{\leq}& 3 C_p \ep^{\frac{p}{2}}+32\sqrt{Q}\sqrt{\pi}\sqrt{\ep}\nn\\
&\leq& 35 C_p \sqrt{\pi}\sqrt{Q} \ep^{\frac{p}{2}}.
\end{eqnarray}
Now we consider first the case $Q\geq 2$. Note
\begin{eqnarray}
s_Q-1&\overset{(\ref{eqvuu1})}{=}&\frac{2+2\frac{Q+1}{Q-1}-\frac{5}{2}-\frac{3}{2}\frac{Q+1}{Q-1}}{\frac{5}{2}+\frac{3}{2}\frac{Q+1}{Q-1}}=\frac{-\frac{1}{2}+\frac{1}{2}\frac{Q+1}{Q-1}}{\frac{5}{2}+\frac{3}{2}\frac{Q+1}{Q-1}}.
\end{eqnarray}
So 
\begin{eqnarray}
\label{peq10}
\frac{s_Q-1}{s_Q}&\overset{(\ref{eqvuu1})}{=}& \frac{\frac{5}{2}+\frac{3}{2}\frac{Q+1}{Q-1}}{2+2\frac{Q+1}{Q-1}}\frac{-\frac{1}{2}+\frac{1}{2}\frac{Q+1}{Q-1}}{\frac{5}{2}+\frac{3}{2}\frac{Q+1}{Q-1}}=\frac{\frac{Q+1}{Q-1}-1}{4\lt(1+\frac{Q+1}{Q-1}\rt)}
=\frac{\lt(\frac{2}{Q-1}\rt)}{4\lt(\frac{2Q}{Q-1}\rt)}=\frac{1}{4Q}.
\end{eqnarray}
Hence 
\begin{eqnarray}
\label{equv25}
\frac{s_Q-1}{\PPI_Q s_Q}&\overset{(\ref{nikc11}),(\ref{peq10})}{=}&\lt(\frac{Q-1}{4Q-2}\rt)\lt(\frac{4Q-1}{2Q}\rt)\frac{1}{4Q}\geq 
\lt(\frac{Q-1}{4Q-2}\rt)\lt(\frac{4Q-2}{2Q}\rt)\frac{1}{4Q}=\frac{Q-1}{8 Q^2}.
\end{eqnarray}

So by Lemma 5.3.1 \cite{astala} and using interpolation of $L^p$ norms (see Section B2, (h) of the Appendix of \cite{evans}) 
and recalling $Q\geq 2$
\begin{eqnarray}
\label{equv30}
\|(\na w_u)_a-(\na w_v)_a\|_{L^2(B_1)}&\leq& \|\mu_{\na w_u}-\mu_{\na w_v}\|_{L^{\frac{P_Q s_Q}{s_Q-1}}(B_1)}\nn\\
&\overset{(\ref{equv3})}{\leq}& 3\lt(\|\mu_{\na w_u}-\mu_{\na w_v}\|_{L^1(B_1)}\rt)^{\frac{s_Q-1}{P_Q s_Q}} \nn\\
&\overset{(\ref{equv24})}{\leq}& 3(35 C_p \sqrt{\pi}\sqrt{Q} \ep^{\frac{p}{2}})^{\frac{s_Q-1}{P_Q s_Q}}\nn\\
&\overset{(\ref{equv25})}{\leq}& 3\pi(35 C_p \sqrt{Q} \ep^{\frac{p}{2}})^{\frac{Q-1}{8 Q^2}}\nn\\
&\leq&6 \pi C_p \sqrt{Q} \ep^{\frac{p}{16 Q^2}}.
\end{eqnarray}
Now in the case $Q<2$ note 
\begin{equation}
\label{peqq10}
\frac{P_Q s_Q}{s_Q-1}\overset{(\ref{eqvuu2}),(\ref{eqvuu1})}{=}3\frac{8}{7}\frac{1}{\frac{8}{7}-1}=24.
\end{equation}
So in the same way as before, using Lemmma 5.3.1 \cite{astala} and interpolation of $L^p$ norms, from 
the second line of (\ref{equv30}) we have 
\begin{eqnarray}
\label{equv15}
\|(D w_u)_a-(D w_v)_a\|_{L^2(B_1)}&\overset{(\ref{equv24}),(\ref{peqq10})}{\leq}& 3(35 C_p \sqrt{\pi}\sqrt{Q} \ep^{\frac{p}{2}})^{\frac{1}{24}}\nn\\
&\leq& 6 \pi C_p\sqrt{Q} \ep^{\frac{p}{48}}\nn\\
&\leq& 6\pi C_p \sqrt{Q} \ep^{\frac{p}{48 Q^2}}.
\end{eqnarray}
Putting (\ref{equv15}) and (\ref{equv30}) together we have 
\begin{equation}
\label{equv17}
\|(D w_u)_a-(D w_v)_a\|_{L^2(B_1)}\leq 12\pi C_p \sqrt{Q}\ep^{\frac{p}{48 Q^2}}.
\end{equation}

Now the Beurling transform $S$ of the anti-conformal part of the gradient of the $L^2$ function 
gives the conformal part of the gradient, see (4.18) Chapter 4 \cite{astala}. So 
$$
S\lt(\frac{\partial}{\partial \bar{z}}(w_u-w_v)\rt)=\frac{\partial}{\partial z}(w_u-w_v).
$$
Since $S$ is an isometry on $L^2(\mathbb{C})$ (using the fact that $w_u$ and $w_v$ are homomorphic outside 
$B_1$ (see (\ref{peq71})) for the last inequality) 
\begin{eqnarray}
\label{nikm2}
\|\frac{\partial w_u}{\partial z}-\frac{\partial w_v}{\partial z}\|_{L^2(\mathbb{C})}&=&
\|S\lt(\frac{\partial w_u}{\partial \bar{z}}\rt)-S\lt(\frac{\partial w_v}{\partial \bar{z}}\rt)\|_{L^2(\mathbb{C})}\nn\\
&\leq&\|\frac{\partial w_u}{\partial \bar{z}}-\frac{\partial w_v}{\partial \bar{z}}\|_{L^2(\mathbb{C})}\nn\\
&=&\frac{1}{\sqrt{2}}\|\lt[D w_u\rt]_a-\lt[D w_v\rt]_a\|_{L^2(B_1)}\nn\\
&\overset{(\ref{equv17})}{\leq}&\frac{12}{\sqrt{2}}\pi C_p \sqrt{Q}\ep^{\frac{p}{48 Q^2}}.\nn
\end{eqnarray}
So 
\begin{equation}
\label{ppeq3}
\|\lt[D w_u\rt]_c-\lt[D w_v\rt]_c\|_{L^2(B_1)}\leq 12 \pi C_p \sqrt{Q}\ep^{\frac{p}{480 Q^2}}.
\end{equation}
Thus 
\begin{eqnarray}
\|\na w_u-\na w_v\|_{L^2(B_1)}&\leq&\|(\na w_u)_a-(\na w_v)_a\|_{L^2(B_1)}
+\|(\na w_u)_s-(\na w_v)_s\|_{L^2(B_1)}
\nn\\
&\overset{(\ref{equv17}),(\ref{ppeq3})}{\leq}& 24\pi C_p \sqrt{Q}\ep^{\frac{p}{48 Q^2}}.  \;\;\;\;\Box\nn
\end{eqnarray}

%
%
%
%

\begin{a1} 
\label{L1.5}
We will show 
\begin{equation}
\label{eaz1}
\|\na w_u\|_{L^2(B_1)}\leq 13\pi Q\text{ and }\|\na w_v\|_{L^2(B_1)}\leq 13\pi Q.
\end{equation} 
Let $\varpi=\min\lt\{3,2+\frac{1}{3(Q-1)}\rt\}$, 
\begin{equation}
\label{nik8}
\|\na w_u\|_{L^{\varpi}(B_1)}\leq 13 Q,\;\|\na w_v\|_{L^{\varpi}(B_1)}\leq 13 Q.
\end{equation}
And 
\begin{equation}
\label{eqz40}
\|w_u-w_v\|_{L^{\infty}(B_{\frac{1}{2}})}\leq 1104 Q^2 \pi C_p \ep^{\frac{p}{120 Q^2}}.
\end{equation}
\end{a1}
\em Proof of Lemma \ref{L1.5}. \rm As before we will take $\kappa=\frac{Q-1}{Q+1}$, so 
\begin{equation}
\label{bbreq51}
\varpi=\min\lt\{3,2+\frac{1-\kappa}{6\kappa}\rt\}\text{ so } \varpi-2=\min\lt\{1,\frac{1-\kappa}{6\kappa}\rt\}.
\end{equation}
Note that 
\begin{equation}
\label{bbreq11}
1-\kappa(1+2(\varpi-2))\overset{(\ref{bbreq51})}{\geq} 1-\kappa\lt(1+\frac{1-\kappa}{3\kappa}\rt)=\frac{2}{3}(1-\kappa) \overset{(\ref{nikc3})}{=}  \frac{4}{3(Q+1)}.
\end{equation}
And note  
\begin{equation}
\label{nik7}
\frac{1-\kappa}{6\kappa}=\frac{1}{6}\lt(\frac{Q+1}{Q-1}\rt)\lt(1-\frac{Q-1}{Q+1}\rt)=\frac{1}{6}\lt(\frac{Q+1}{Q-1}\rt)\lt(\frac{2}{Q+1}\rt)=\frac{1}{3(Q-1)}
\end{equation}
and thus 
\begin{equation}
\label{eqaa1}
1\overset{(\ref{bbreq51})}{\geq} \varpi-2\overset{(\ref{nik7}),(\ref{bbreq51})}{=}\min\lt\{1,\frac{1}{3(Q-1)}\rt\}\geq \frac{1}{3Q}.
\end{equation}
So from Lemma \ref{LA2} (\ref{eqc7}) (and recalling (\ref{bbreq51}))
\begin{eqnarray}
\label{nikc56}
\lt(\int_{B_1} \lt|\na w_u\rt|^\varpi \rt)^{\frac{1}{\varpi}}
&\overset{(\ref{eqc7}),(\ref{bbreq11})}{\leq}& 1+\frac{3}{2}(Q+1)(1+3(\varpi-2))\\
&\overset{(\ref{eqaa1})}{\leq}& 
1+6(Q+1)\nn\\
&\leq& 13Q.
\end{eqnarray}
In the same way $\|D w_v\|_{L^{\varpi}(B_1)}\leq 13 Q$, thus (\ref{nik8}) is established. 
By Holder if we let $r=\frac{\varpi}{2}$ and $r'>0$ be such that $\frac{1}{r}+\frac{1}{r'}=1$ and so 
\begin{eqnarray}
\int_{B_1} \lt|\na w_u\rt|^2 dz&\leq& \lt(\int_{B_1} \lt|\na w_u\rt|^{2r} dz\rt)^{\frac{1}{r}}\lt(\int_{B_1} 1 dz\rt)^{\frac{1}{r'}}\nn\\
&\leq& \pi \lt( \int_{B_1} \lt|\na w_u\rt|^\varpi dz\rt)^{\frac{1}{r}}\nn\\
&\overset{(\ref{nik8})}{\leq}& \pi \lt(13 Q\rt)^{\frac{\varpi}{r}}\nn\\
&\leq& 13^2 \pi Q^2.\nn
\end{eqnarray}
So $\|\na w_u\|_{L^2(B_1)}\leq 13\pi Q$ and in the same way 
$\|\na w_v\|_{L^2(B_1)}\leq 13\pi Q$. So (\ref{eaz1}) is established.  

Let 
\begin{equation}
\label{eqc66}
r=1+\frac{\varpi}{2}=\frac{2+\varpi}{2}.
\end{equation}
Since by (\ref{bbreq51}), $2<\varpi\leq 3$, so 
$r\in (2,\varpi)$ and thus $\frac{1}{\varpi}<\frac{1}{r}<\frac{1}{2}$ and thus 
there exists $\theta\in (0,1)$ such that 
\begin{equation}
\label{bbreq27}
\frac{1}{r}=\frac{\theta}{2}+\frac{1-\theta}{\varpi}. 
\end{equation}
By interpolation of $L^p$ norms we know 
\begin{eqnarray}
\label{nik201}
\|\na w_u-\na w_v\|_{L^r(B_1)}&\leq& \|\na w_u-\na w_v\|^{\theta}_{L^2(B_1)}
\|\na w_u-\na w_v\|_{L^{\varpi}(B_1)}^{1-\theta}\nn\\
&\overset{(\ref{equa127})}{\leq}& \lt(24 \pi C_p \sqrt{Q} \ep^{\frac{p}{48 Q^2}}\rt)^{\theta}\lt(\|\na w_v\|_{L^{\varpi}(B_1)}+\|\na w_u\|_{L^{\varpi}(B_1)}\rt)^{1-\theta}\nn\\
&\overset{(\ref{eaz1})}{\leq}& \lt( 24\pi C_p \sqrt{Q} \ep^{\frac{p}{48 Q^2}}  \rt)^{\theta}
\lt(26\pi Q\rt)^{1-\theta}.
\end{eqnarray}
Now since $r=\frac{2+\varpi}{2}$ so $\frac{2}{2+\varpi}-\frac{1}{\varpi}\overset{(\ref{bbreq27})}{=}\theta(\frac{1}{2}-\frac{1}{\varpi})$. So  
\begin{eqnarray}
\label{ffefq1}
\theta\lt(\frac{\varpi-2}{2\varpi}\rt)&=&\theta\lt(\frac{1}{2}-\frac{1}{\varpi}\rt)\nn\\
&=&\frac{2}{2+\varpi}-\frac{1}{\varpi}\nn\\
&=& \frac{\varpi-2}{\varpi(2+\varpi)}.
\end{eqnarray}
So again since by (\ref{bbreq51}) $2< \varpi\leq 3$, thus 
$\frac{2}{5}\leq \theta\overset{(\ref{ffefq1})}{=}\frac{2}{(2+\varpi)}< \frac{1}{2}$. Thus 
\begin{eqnarray}
\label{eqz1}
\|\na w_u-\na w_v\|_{L^r(B_1)}&\overset{(\ref{nik201})}{\leq}& 24^{\theta}\pi (26)^{1-\theta} C_p Q^{\theta} Q^{1-\theta} \ep^{\frac{p}{120 Q^2}}\nn\\
&\leq& 26 Q\pi C_p  \ep^{\frac{p}{120 Q^2}}.
\end{eqnarray}

Now from the proof of Lemma 4.28 of \cite{adams} letting $Q_r(x)$ denote the square of side length $r$ centred on $x$, we have 
that 
\begin{eqnarray}
\label{eqz15}
\lt|(w_u-w_v)(x)-\frac{1}{2}\int_{Q_{\frac{1}{\sqrt{2}}}(0)} (w_u-w_v)(z) dz \rt|\leq 
K \lt(\frac{1}{\sqrt{2}}\rt)^{1-\frac{2}{r}} \|\na w_u-\na w_v\|_{L^r(B_1)}
\end{eqnarray}
where 
\begin{eqnarray}
\label{eqz23}
K&=&\sqrt{2}\int_{0}^1 t^{-\frac{2}{r}} dt\overset{(\ref{eqc66})}{=}\sqrt{2}\int_{0}^1 t^{-\frac{4}{2+\varpi}} dt\nn\\
&=&\sqrt{2}  \lt(\frac{2+\varpi}{\varpi-2}\rt) 
\int_{0}^1 \frac{d}{dt}\lt(t^{\frac{\varpi-2}{2+\varpi}}\rt) dt=\sqrt{2}\lt(\frac{2+\varpi}{\varpi-2}\rt)  \overset{(\ref{eqaa1}), (\ref{bbreq51})}{\leq}15 \sqrt{2}Q.
\end{eqnarray}

So by (\ref{eqz1}), (\ref{eqz15}), (\ref{eqz23}) we have that 
\begin{eqnarray}
\lt|(w_u-w_v)(x)-\frac{1}{2}\int_{Q_{\frac{1}{\sqrt{2}}}(0)} (w_u-w_v)(z) dz\rt|&\leq& 15\sqrt{2}Q\times 26 \pi Q C_p  \ep^{\frac{p}{120 Q^2}}\nn\\
&\leq& 390 \sqrt{2} \pi Q^2 C_p  \ep^{\frac{p}{120 Q^2}}.\nn
\end{eqnarray}

So 
$$
\lt|(w_u-w_v)(x)-(w_u-w_v)(y)\rt| \leq 780 \sqrt{2} \pi Q^2 C_p  \ep^{\frac{p}{120 Q^2}}    \text{ for 
any }x,y\in Q_{\frac{1}{\sqrt{2}}}(0).
$$
This establishes (\ref{eqz40}). $\Box$

%
%
%
%

\begin{a1} 
\label{L2}
Given $Q$-quasiregular mappings $u$, $v$ with the property that 
\begin{equation}
\label{eq17}
\int_{B_1} \lt|S(\na u)-S(\na v)\rt|^2 dz\leq \ep
\end{equation}
 then letting $w_u,\phi_u$ denote the 
Stoilow decomposition of $u$ and $w_v,\phi_v$ denote the Stoilow decomposition of $v$. We will show that 
\begin{equation}
\label{zueq5}
\int_{w_u(B_{\frac{1}{2}})} \lt|\lt|\phi_u'(y)\rt|^2-\lt|\phi_v'(y)\rt|^2\rt| dy \leq 2.3\times 10^9 \pi^5 Q^4 \mu^{-3} C_p \ep^{\frac{p}{120 Q^2}}.
\end{equation}
\end{a1}

\em Proof of Lemma \ref{L2}. \rm Note $\na u(z)=\na \phi_u(w_u(z))\na w_u(z)$, so 
\begin{equation}
\label{mfeq1}
\na u(z)^T \na u(z)=\na w_u(z)^T\na \phi_u(w_u(z))^T \na\phi_u(w_u(z)) \na w_u(z)=\lt|\na \phi_u(w_u(z))\rt|^2 \na w_u(z)^T \na w_u(z).
\end{equation}
We know 
\begin{equation}
\label{kleq1}
\phi_v'(z)=\mathrm{Re}\lt(\phi_v(z)\rt)_x+i\mathrm{Im}\lt(\phi_v(z)\rt)_x=
\mathrm{Im}\lt(\phi_v(z)\rt)_y-i \mathrm{Re}\lt(\phi_v(z)\rt)_y.
\end{equation}
So to simplify notation let 
\begin{equation}
\label{bbreq81}
\lm(z)=\lt|\na \phi_u(z)\rt|^2=2\lt|\phi'_u(z)\rt|^2\text{ and }\varrho(z)\overset{(\ref{kleq1})}{=}\lt|\na \phi_v(z)\rt|^2=2\lt|\phi'_v(z)\rt|^2. 
\end{equation}
Thus from (\ref{mfeq1}), (\ref{bbreq81})
\begin{equation}
\label{bbreq68}
\na u(z)^T \na u(z)=\lm(w_u(z))\na w_u(z)^T \na w_u(z)\text{ and }\na v(z)^T \na v(z)=\varrho(w_v(z))\na w_v(z)^T \na w_v(z).
\end{equation}

Note 
$$
\lt|S(Du)\rt|\overset{(\ref{nikb1})}{\leq} 2\|S(Du)\|=2\|Du\|<2\lt|Du\rt|,
$$
so 
\begin{equation}
\label{akeq1}
\|S(Du)\|_{L^2(B_1)}\leq 2\text{ and }\|S(Dv)\|_{L^2(B_1)}\overset{(\ref{peq501})}{\leq} 3.
\end{equation}
Thus 
\begin{eqnarray}
\label{kleq11}
&~&\int_{B_1} \lt|S(Du)^2-S(Dv)^2\rt| dz\nn\\
&~&\qd\qd\leq 
\int_{B_1} \lt|S(Du)\lt(S(Du)-S(Dv)\rt)\rt|+\lt|\lt(S(Du)-S(Dv)\rt)S(Dv)\rt| dz\nn\\
&~&\qd\qd\leq \|S(Du)\|_{L^2(B_1)}\|S(Du)-S(Dv)\|_{L^2(B_1)}+ \|S(Dv)\|_{L^2(B_1)}\|S(Du)-S(Dv)\|_{L^2(B_1)}\nn\\
&~&\qd\qd \overset{(\ref{akeq1})}{\leq} 5\sqrt{\ep}.
\end{eqnarray}
Recall constant $\mu=\mucon$ and $\gamma=\gammacon$. Since $S(\na u)^2=\na u^T \na u$ and $S(\na v)^2=\na v^T \na v$
\begin{eqnarray}
\label{eq20.4}
10\sqrt{\ep}&\overset{(\ref{kleq11})}{\geq}& \int_{B_{\gamma}} \lt|Tr(\na u^T \na u)-Tr(\na v^T \na v)\rt| dz\nn\\
&\overset{(\ref{bbreq68})}{=}&\int_{B_{\gamma}} \lt| \lm(w_u) Tr(\na w_u^T \na w_u)-\varrho(w_v)Tr(\na w_v^T \na w_v)\rt| dz\nn\\
&\geq& \int_{B_{\gamma}} \lt| \lm(w_u) Tr(\na w_u^T \na w_u)-\varrho(w_u)Tr(\na w_u^T \na w_u)\rt| dz\nn\\
&~&-\int_{B_{\gamma}} \lt| \varrho(w_u) Tr(\na w_u^T \na w_u)-\varrho(w_v)Tr(\na w_v^T \na w_v)\rt| dz.
\end{eqnarray}
Now note 
\begin{equation}
\label{ezq6}
\sup\lt\{\lt|\varrho(w_u(z))\rt|:x\in B_{\gamma}\rt\}\overset{(\ref{ezq1})}{\leq}  \sup\lt\{\lt|\varrho(y)\rt|:y\in B_{\mu}(w_u(0))\rt\}\overset{(\ref{bbreq81}),(\ref{eq3})}{\leq} 
2\times \lt(\frac{16 \pi}{\mu}\rt)^2=\frac{512 \pi^2}{\mu^2}.
\end{equation}
Now from (\ref{kleq1}) as 
$$
\phi''_v=\mathrm{Re}(\phi_v)_{xx}+i \mathrm{Im}(\phi_v)_{xx}= \mathrm{Im}(\phi_v)_{xy}-i \mathrm{Re}(\phi_v)_{xy}=-\mathrm{Re}(\phi_v)_{yy}-i \mathrm{Im}(\phi_v)_{yy}.
$$
Thus we have 
\begin{equation}
\label{pw2}
\lt|D^2\phi_v(z)\rt|\leq 4 \lt|\phi_v''(z)\rt|.
\end{equation}
And note for $k=1,2$, 
\begin{eqnarray}
\lt|\varrho_{,k}(z)\rt|&\overset{(\ref{bbreq81})}{\leq}& 2\lt|D^2 \phi_v(z)\rt|\lt|D\phi_v(z)\rt|\nn\\
&\overset{(\ref{pw2}),(\ref{bbreq81})}{\leq}& 16\lt|\phi_v''(z)\rt|\lt|\phi_v'(z)\rt|\nn\\
&\overset{(\ref{eq3}),(\ref{eq5})}{\leq}& \frac{12288\pi^2}{\mu^3}\text{ for any }z\in B_{\mu}(w_u(0)). 
\end{eqnarray}
Thus 
\begin{equation}
\label{nikg3}
\lt|D\varrho(z)\rt|\leq 
{\frac{24576\pi^2}{\mu^3}}\text{ for any }z\in B_{\mu}(w_u(0)).
\end{equation}
Hence 
\begin{eqnarray}
\label{ecx31}
\sup\lt\{\lt|\varrho(w_u(z))-\varrho(w_v(z))\rt|:x\in B_{\gamma}\rt\}&\overset{(\ref{ezq1}),(\ref{nikg3})}{\leq}&
\frac{24576\pi^2}{\mu^3}\sup\lt\{\lt|w_u(y)-w_v(y)\rt|:y\in B_{\gamma}\rt\}  \nn\\
&\overset{(\ref{eqz40})}{\leq}& 2.72 \times 10^7 Q^2 C_p \frac{\pi^3}{\mu^3}\ep^{\frac{p}{120 Q^2}}.
\end{eqnarray}
Thus
\begin{eqnarray}
\label{ezq10}
&~&\int_{B_{\gamma}} \lt| \varrho(w_u) Tr(\na w_u^T \na w_u)-\varrho(w_v)Tr(\na w_v^T \na w_v)\rt| dz\nn\\ 
&~&\leq \int_{B_{\gamma}} \lt| \varrho(w_u) Tr(\na w_u^T \na w_u)-\varrho(w_u)Tr(\na w_v^T \na w_v)\rt| dz\nn\\ 
&~&\qd\qd\qd + \lt|(\varrho(w_u)-\varrho(w_v))Tr(\na w_v^T \na w_v)\rt| dz\nn\\
&~&\overset{(\ref{ezq6}),(\ref{ecx31})}{\leq} \frac{1024\pi^2}{\mu^2}\int_{B_{\gamma}} \lt|\na w_u^T \na w_u-\na w_v^T \na w_v\rt| dz+
2.72 \times 10^7 C_p Q^2\frac{\pi^3}{\mu^3} \ep^{\frac{p}{120 Q^2}}\int_{B_{\frac{1}{2}}} \lt|\na w_v\rt|^2 dz.\nn\\
\end{eqnarray}
Now 
\begin{eqnarray}
\label{ezq11}
\int_{B_1} \lt|\na w_u^T \na w_u-\na w_v^T \na w_v\rt| dz&\leq& 
\int_{B_1} \lt|\na w_u^T (\na w_u-\na w_v)\rt| dz+
\int_{B_1} \lt|(\na w_u^T-\na w_v^T)\na w_v\rt| dz\nn\\
&\overset{(\ref{eaz1})}{\leq}& 26 \pi Q\|\na w_u-\na w_v\|_{L^2(B_1)}\nn\\
&\overset{(\ref{equa127})}{\leq}& 624 \pi^2 C_p Q^{2} \ep^{\frac{p}{48 Q^2}}. 
\end{eqnarray}
So applying (\ref{ezq11}) and (\ref{eaz1}) to (\ref{ezq10})
\begin{eqnarray}
&~&
\int_{B_{\gamma}} \lt|\varrho(w_u)Tr(\na w_u^T \na w_u)- \varrho(w_v)Tr(\na w_v^T \na w_v)\rt| dz \nn\\
&~&\qd\qd\qd\qd\overset{(\ref{eaz1}),(\ref{ezq10}),(\ref{ezq11})}{\leq} 
\frac{1024 \pi^2}{\mu^2}\times 624 \pi^2 C_p Q^{2} \ep^{\frac{p}{48 Q^2}}
+2.72 \times 10^7 C_p Q^2 \frac{\pi^3}{\mu^3} \ep^{\frac{p}{120 Q^2}}\times 13^2 \pi^2 Q^2\nn\\
&~&\qd\qd\qd\qd\leq 4.598\times 10^9 \mu^{-3}\pi^5 Q^4 C_p \ep^{\frac{p}{120 Q^2}}.\nn
\end{eqnarray}
Putting this together with (\ref{eq20.4}) we have that 
\begin{eqnarray}
 4.6 \times 10^9 \mu^{-3}\pi^5 Q^4 C_p \ep^{\frac{p}{120 Q^2}}&\geq& \int_{B_{\frac{1}{2}}} \lt| (\lm(w_u)-\varrho(w_u))Tr(\na w_u^T \na w_u)\rt| dz\nn\\
&\geq& \int_{B_{\frac{1}{2}}} \lt|\lm(w_u)-\varrho(w_u)\rt|\mathrm{det}(\na w_u) dz\nn\\
&\overset{(\ref{bbreq81})}{=}&2\int_{w_u\lt(B_{\frac{1}{2}}\rt)} \lt|\lt|\phi_u'(y)\rt|^2-\lt|\phi_v'(y)\rt|^2\rt| dy\;\;\;\;\;\;\Box\nn
\end{eqnarray}

%
%
%
%

\begin{a1} 
\label{L3}
Let
\begin{equation}
\label{peq51}
\mu=\mucon. 
\end{equation}
Recall that $B_{2\mu}(w_u(0))\subset w_u(B_{\frac{1}{2}}(0))$. Fix 
constant 
\begin{equation}
\label{breq80}
h_0=\frac{\mu^2}{96\pi}\lt(\frac{\mucsta}{\mucstc}\rt)^{\frac{1}{2p}}
\end{equation}
We can 
find $x_0\in B_{\frac{\mu}{2}}(w_u(0))$ such that  
\begin{equation}
\label{breq60}
\inf\lt\{\lt|\phi_u'(y)\rt|:y\in B_{h_0}(x_0)\rt\}\geq \frac{1}{2}
\lt(\frac{\mucsta}{\mucstc}\rt)^{\frac{1}{2p}}.
\end{equation}
\end{a1}

\em Proof of Lemma \ref{L3}. \rm 

Note 
\begin{eqnarray}
\label{nikq1}
C_p&\geq& \int_{B_1} \det(Du (z))^{-p} dz\nn\\
&=& \int_{B_1} 
\det(D\phi_u(w_u(z)))^{-p}\det(Dw_u(z))^{-p} dz\nn\\
&\overset{(\ref{eqz11})}{\geq}&\int_{w_u^{-1}(B_{\mu}(w_u(0)))} \det(D\phi_u(w_u(z)))^{-p}\det(Dw_u(z)) \det(Dw_u(z))^{-p-1} dz\nn\\
&=&\int_{w_u^{-1}(B_{\mu}(w_u(0)))} \det(D\phi_u(w_u(z)))^{-p}\det(Dw_u(w_u^{-1}(w_u(z))))^{-p-1} \det(Dw_u(z)) dz\nn\\
&=&\int_{B_{\mu}(w_u(0))} \det(D\phi_u(y))^{-p}\det(Dw_u(w_u^{-1}(y)))^{-p-1} dy.
\end{eqnarray}
Let $\varsigma>4Q$ be some constant we decide on later
\begin{equation}
\label{mareq2}
D_{\varsigma}=\lt\{z\in B_1: \det(Dw_u(z))>\varsigma\rt\}. 
\end{equation}
Thus by Theorem 13.1.4 \cite{astala}
$$
Q \pi\lt(\frac{\lt|D_{\varsigma}\rt|}{\pi}\rt)^{\frac{1}{Q}}\geq \int_{D_{\varsigma}} \det(Dw_u(z)) dz\geq 
\varsigma \lt|D_{\varsigma}\rt|.
$$
So $\frac{\lt|D_{\varsigma}\rt|}{\pi}\geq \lt(\frac{\varsigma}{Q \pi}\rt)^Q\lt|D_{\varsigma}\rt|^Q$ and 
thus $\lt|D_{\varsigma}\rt|^{Q-1}\leq \frac{Q^Q \pi^{Q-1}}{\varsigma^Q}$. Hence as $\varsigma>4Q$
\begin{equation}
\label{mareq2.7}
\lt|D_{\varsigma}\rt|\leq \pi\lt(\frac{Q}{\varsigma}\rt)^{\frac{Q}{Q-1}}\leq \pi\lt(\frac{Q}{\varsigma}\rt).
\end{equation}
In particular 
\begin{equation}
\label{peq120}
\lt|D_{\varsigma}\rt|<1.
\end{equation}
Now let $\varphi=\min\lt\{\frac{3}{2},1+\frac{1}{6(Q-1)}\rt\}$. Note $\varphi\geq 1+\frac{1}{6Q}$. 
\begin{equation}
\label{nikzz5}
\frac{\varphi-1}{\varphi}\geq \frac{1}{6Q\varphi}\geq \frac{1}{9Q}.
\end{equation}
Now note $2\varphi=\varpi$ where 
$\varpi$ is the constant from from the statement of Lemma \ref{L1.5}. Note 
\begin{eqnarray}
\label{nik15}
\int_{D_{\varsigma}} \det(Dw_u) dz&\leq& \int_{D_{\varsigma}} \lt| Dw_u\rt|^2 dz\nn\\
&\leq& \lt( \int_{D_{\varsigma}} \lt| D w_u\rt|^{2\varphi} dz\rt)^{\frac{1}{\varphi}}\lt|D_{\varsigma}\rt|^{\frac{1}{\varphi'}}\nn\\
&=&\lt(\lt( \int_{D_{\varsigma}} \lt| D w_u\rt|^{\varpi} dz\rt)^{\frac{1}{\varpi}}\rt)^2
\lt|D_{\varsigma}\rt|^{\frac{1}{\varphi'}}\nn\\
&\overset{(\ref{nik8})}{\leq}& (13 Q)^{2}\lt|D_{\varsigma}\rt|^{\frac{\varphi-1}{\varphi}}\nn\\
&\overset{(\ref{peq120}),(\ref{nikzz5})}{\leq}& (13 Q)^{2}\lt|D_{\varsigma}\rt|^{\frac{1}{9 Q}}\nn\\
&\overset{(\ref{mareq2.7})}{\leq}& (13 Q)^2 \pi \lt(\frac{Q}{\varpi}\rt)^{\frac{1}{9 Q}}.
\end{eqnarray}
Now let 
\begin{equation}
\label{nik16}
\varsigma=(1352)^{9 Q} Q^{(18Q+1)}\mu^{-18 Q},
\end{equation}
so 
\begin{equation}
\label{nik17}
\varsigma^{\frac{1}{9Q}}=1352 Q^{\frac{1}{9Q}+2}\mu^{-2},
\end{equation}
thus 
\begin{equation}
\label{nik18}
\mu^2=1352 Q^2 \lt(\frac{Q}{\varsigma}\rt)^{\frac{1}{9Q}},
\end{equation}
hence 
\begin{equation}
\label{nik19}
\pi\mu^2 8^{-1}=(13 Q)^2 \pi\lt(\frac{Q}{\varsigma}\rt)^{\frac{1}{9Q}}.
\end{equation}
So note by (\ref{nik15}) we have that 
\begin{eqnarray}
\label{nik20}
\lt|B_{\frac{\mu}{2}}(w_u(0))\backslash w_u(D_{\varsigma})\rt|&\overset{(\ref{nik15})}{\geq}& \pi\frac{\mu^2}{4}-(13 Q)^2 \pi \lt(\frac{Q}{\varsigma}\rt)^{\frac{1}{9Q}}\nn\\
&\overset{(\ref{nik19})}{\geq}&\pi\frac{\mu^2}{8}.
\end{eqnarray}
Thus 
\begin{eqnarray}
\label{nik31}
C_p&\overset{(\ref{nikq1})}{\geq}&\int_{B_{\frac{\mu}{2}}(w_u(0))\backslash w_u(D_{\varsigma})}  \det(D\phi_u(y))^{-p}\det(Dw_u(w_u^{-1}(y)))^{-p-1} dy\nn\\
&\geq&\varsigma^{-p-1} \int_{B_{\frac{\mu}{2}}(w_u(0))\backslash w_u(D_{\varsigma})}  \det(D\phi_u(y))^{-p} dy\\
&\overset{(\ref{nik20})}{\geq}&\inf\lt\{\det(D\phi_u(y))^{-p}:y\in B_{\frac{\mu}{2}}(w_u(0))\backslash w_u(D_{\varsigma})\rt\} \pi \varsigma^{-p-1}\frac{\mu^2}{8}.
\end{eqnarray}
So there must exist $\zeta_0\in B_{\frac{\mu}{2}}(w_u(0))\backslash w_u(D_{\varsigma})$ such that 
$\det(D\phi_u(\zeta_0))^{-p}\leq \frac{9 C_p \varsigma^{p+1}}{\mu^2}$. Note 
\begin{equation}
\label{peq20}
\varsigma^2\overset{(\ref{nik16})}{\leq} (1352)^{18 Q} Q^{40 Q} \mu^{-36 Q},
\end{equation} 
thus (recalling $p\in (0,1)$)
\begin{equation}
\det(D\phi_u(\zeta_0))\geq\frac{\mu^{\frac{2}{p}}}{9^{\frac{1}{p}}C_p^{\frac{1}{p}}\varsigma^{\frac{p+1}{p}}}\geq \frac{\mu^{\frac{2}{p}}}{ 9^{\frac{1}{p}}C_p^{\frac{1}{p}}\varsigma^{\frac{2}{p}}}\overset{(\ref{peq20})}{\geq} 
\lt(\frac{\mucsta}{\mucstc}\rt)^{\frac{1}{p}}.
\end{equation}
So if $y\in B_{h_0}(\zeta_0)$ then 
\begin{eqnarray}
\lt|\phi_u'(y)-\phi_u'(\zeta_0)\rt|&\leq& \int_{\lt[y,x_0\rt]} \lt|\phi_u''(z)\rt| dH^1 z\nn\\
&\overset{(\ref{eq5})}{\leq} &h_0\frac{48\pi}{\mu^2}. 
\end{eqnarray}
Hence as $h_0\overset{(\ref{breq80})}{=}\frac{\mu^2}{96 \pi} \lt(\frac{\mucsta}{\mucstc}\rt)^{\frac{1}{2p}}$
\begin{eqnarray}
\lt|\phi_u'(y)\rt|&\geq&  \lt(\frac{\mucsta}{\mucstc}\rt)^{\frac{1}{2p}}
-\frac{48 h_0\pi}{\mu^2}\nn\\
&\geq& \frac{1}{2} \lt(\frac{\mucsta}{\mucstc}\rt)^{\frac{1}{2p}}   \text{ for any }y\in B_{h_0}(x_0).\;\;\;\Box \nn
\end{eqnarray}

%
%
%
%

\begin{a1} 
\label{L4}

We will show there exists $\zeta\in \mathbb{C}$ such that 
\begin{equation}
\label{eq60}
\sup\lt\{\lt|\phi_u'(z)-\zeta \phi_v'(z)\rt|:z\in B_{\frac{\mu}{2}}(w_u(0))\rt\}\leq 
C_p \CI_1 \ep^{\frac{p^3}{4\times 10^{7} Q^5 \log(10 C_p Q)}}.
\end{equation}
\end{a1}

\em Proof of Lemma \ref{L4}. \rm Let $h_0$ be the constant defined by (\ref{breq80}) of Lemma \ref{L3} and let $x_0\in B_{\frac{\mu}{2}}(w_u(0))$ be the point from Lemma \ref{L3} that satisfies (\ref{breq60}). 

Note since $x_0\in B_{\frac{\mu}{2}}(w_u(0))$ and $h_0\overset{(\ref{breq80})}{\leq} \frac{\mu}{2}$ thus 
\begin{equation}
\label{hjeq1}
B_{h_0}(x_0)\subset B_{\mu}(w_u(0))\overset{(\ref{eqz11})}{\subset} w_u\lt(B_{\frac{1}{2}}\rt).
\end{equation}

Now 
\begin{eqnarray}
2.3\times 10^9 \pi^5 Q^4 \mu^{-3} C_p \ep^{\frac{p}{120 Q^2}}&\overset{(\ref{zueq5})}{\geq}& \int_{B_{h_0}(x_0)} 
\lt|\lt|\phi_u'(y)\rt|^2-\lt|\phi_v'(y)\rt|^2\rt| dy\nn\\
&=&\int_{B_{h_0}(x_0)} \lt|\lt|\phi_u'(y)\rt|-\lt|\phi_v'(y)\rt|\rt|\lt|\lt|\phi_u'(y)\rt|+\lt|\phi_v'(y)\rt|\rt|  dy\nn\\
&\overset{(\ref{breq60})}{\geq}&  \frac{1}{2}\lt(\frac{\mucsta}{\mucstc}\rt)^{\frac{1}{2p}}\int_{B_{h_0}(x_0)} \lt|\lt|\phi_u'(y)\rt|-\lt|\phi_v'(y)\rt|\rt| dy.\nn
\end{eqnarray}
Thus 
\begin{equation}
\label{fefq1}
\int_{B_{h_0}(x_0)} 
\lt|\lt|\phi_u'(y)\rt|-\lt|\phi_v'(y)\rt|\rt| dy\leq \lt(\frac{\mucstc}{\mucsta} \rt)^{\frac{1}{2p}} 4.6\times 10^9 \pi^5 Q^4 \mu^{-3} C_p \ep^{\frac{p}{120 Q^2}}. 
\end{equation}
By Cauchy's theorem we can find an analytic function $\psi$ such that 
\begin{equation}
\label{peq170}
\psi'(z)=\frac{\phi_v'(z)}{\phi_u'(z)}\text{ for }z\in B_{h_0}(x_0). 
\end{equation}
So 
\begin{eqnarray}
\int_{B_{h_0}(x_0)} \lt|1-\lt|\psi'(z)\rt|\rt|^2 dz&=&\int_{B_{h_0}(x_0)} \lt|\phi_u'(z)\rt|^{-2}
\lt|\lt|\phi_u'(z)\rt|- \lt|\phi_v'(z)\rt|\rt|^2 dz\nn\\
&\overset{(\ref{breq60}),(\ref{eq3})}{\leq}& \frac{128\pi}{\mu}   
\lt(\frac{\mucsta}{\mucstc}\rt)^{-\frac{1}{p}}  \int_{B_{h_0}(x_0)} \lt|\lt|\phi_u'(z)\rt|- \lt|\phi_v'(z)\rt|\rt| dz\nn\\
&\overset{(\ref{fefq1})}{\leq}&  6\times 10^{11} \lt(\frac{\mucstc}{\mucsta}\rt)^{\frac{3}{2p}}  \pi^6 Q^{4} \mu^{-4} C_p \ep^{\frac{p}{120 Q^2}}.
\end{eqnarray}
Now since $\lt[\psi'(z)\rt]_M\in CO_{+}(2)$, $\sqrt{2}\lt|1-\lt|\psi'(z)\rt|\rt|\overset{(\ref{olleq6})}{=}\mathrm{dist}(D\psi(z),SO(2))$. So 
\begin{eqnarray}
\label{breq120}
\Xint{-}_{B_{h_0}(x_0)} \mathrm{dist}^2\lt(D\psi(z),SO(2)\rt) dz&\leq& 
 \frac{12\times 10^{11}}{h_0^2}  \lt(\frac{\mucstc}{\mucsta}\rt)^{\frac{3}{p}} \pi^5 Q^{4} \mu^{-4} C_p \ep^{\frac{p}{120Q^2}}\nn\\
&\overset{(\ref{breq80})}{\leq}& (96\pi)^2\lt(\frac{9 C_p (1352)^{18Q} Q^{40Q}}{\mu^{38Q}}\rt)^{\frac{4}{p}}
12\times 10^{11} \pi^5 Q^4 \mu^{-8} C_p \ep^{\frac{p}{120Q^2}}\nn\\
&\leq& \lt(\sqrt{96\pi}\times 12^{\frac{1}{4}} \times 10^{\frac{11}{4}} \times 9 \times (1352)^{18} \pi^2\rt)^{\frac{4Q}{p}}
\lt(\frac{Q}{\mu}\rt)^{\frac{164Q}{p}} C_p^{\frac{5}{p}} \ep^{\frac{p}{120Q^2}}\nn\\
&\leq& (\cfc)^{\frac{4Q}{p}}
\lt(\frac{Q}{\mu}\rt)^{\frac{164 Q}{p}} C_p^{\frac{5}{p}} \ep^{\frac{p}{120Q^2}}.
\end{eqnarray}
Let $\zeta(z)=\psi(x_0+h_0 z)h_0^{-1}$. Thus 
$$
\int_{B_1} \mathrm{dist}^2\lt(D\zeta(z),SO(2)\rt) dz\leq (\cfc)^{\frac{4Q}{p}}
\lt(\frac{Q}{\mu}\rt)^{\frac{164Q}{p}} C_p^{\frac{5}{p}} \ep^{\frac{p}{120Q^2}}.
$$
So in particular $\|D\zeta \|_{L^2(B_1)}\leq 2$. Thus by applying Proposition \ref{PP1} we have that 
there exists $R\in SO(2)$ such that 
$$
\int_{B_{\frac{1}{4}}} \lt|D\zeta(z) -R\rt|^2 dz\leq 9 \times (\cfc)^{\frac{Q}{p}}
\lt(\frac{Q}{\mu}\rt)^{\frac{41Q}{p}} C_p^{\frac{2}{p}}\ep^{\frac{p}{480 Q^2}}.
$$
By rescaling we obtain that there exists $R$ such that 
$$
\Xint{-}_{B_{\frac{h_0}{4}}} \lt|D\psi(z)-R\rt|^2 dz\leq 9 \times (\cfc)^{\frac{Q}{p}}
\lt(\frac{Q}{\mu}\rt)^{\frac{41Q}{p}} C_p^{\frac{2}{p}}\ep^{\frac{p}{480 Q^2}}.
$$
Thus Holder's inequality 
\begin{equation}
\label{nikzz33}
\Xint{-}_{B_{\frac{h_0}{4}}}    \lt|D\psi(z)-R\rt| dz\leq  3\times (\cfc)^{\frac{Q}{2p}}
\lt(\frac{Q}{\mu}\rt)^{\frac{21Q}{p}} C_p^{\frac{1}{p}}\ep^{\frac{p}{960 Q^2}}. 
\end{equation}
Returning to complex notation for some $\zeta_1\in \mathbb{C}\cap \lt\{z:\lt|z\rt|=1\rt\}$ we have 
\begin{eqnarray}
\label{nikzz11}
\int_{B_{\frac{h_0}{4}}} \lt|\psi'(z)-\zeta_1\rt| dz&\leq&  48 h_0^2 \times (\cfc)^{\frac{Q}{2p}}
\lt(\frac{Q}{\mu}\rt)^{\frac{21Q}{p}} C_p^{\frac{1}{p}}\ep^{\frac{p}{960 Q^2}}. 
\end{eqnarray}
Now by the Co-area formula we know 
\begin{equation}
\label{nikw3}
\int_{\frac{h_0}{8}}^{\frac{h_0}{4}} \int_{\partial B_s(x_0)} \lt|\psi'(z)-\zeta_1\rt| dH^1 x \;ds\leq \int_{B_{\frac{h_0}{4}}(x_0)} \lt|\psi'(z)-\zeta_1\rt| dz
\end{equation}
So we must be able to find 
\begin{equation}
\label{oopeq70}
q\in \lt(\frac{h_0}{8},\frac{h_0}{4}\rt) 
\end{equation}
such that 
\begin{eqnarray}
\label{peq25}
\int_{\partial B_q(x_0)} \lt|\psi'(z)-\zeta_1\rt| dH^1 z&\leq& \frac{8}{h_0} \int_{B_{\frac{h_0}{4}}(x_0)} \lt|\psi'(z)-\zeta_1\rt| dz\nn\\
&\overset{(\ref{nikzz11}),(\ref{nikw3})}{\leq}& 384 h_0 \times (\cfc)^{\frac{Q}{2p}}
\lt(\frac{Q}{\mu}\rt)^{\frac{21Q}{p}} C_p^{\frac{1}{p}}\ep^{\frac{p}{960 Q^2}}. 
\end{eqnarray}
So
\begin{eqnarray}
\label{jejq1}
\int_{\partial B_q(x_0)} \lt|\phi_v'(z)-\zeta_1 \phi_u'(z)\rt| dH^1 z&\overset{(\ref{peq170})}{=}&
 \int_{\partial B_q(x_0)} \lt|(\psi'(z)-\zeta_1)\phi_u'(z)\rt| dH^1 z\nn\\
&\overset{(\ref{eq3})}{\leq}&\frac{16 \pi}{\mu} \int_{\partial B_q(x_0)}  \lt|\psi'(z)-\zeta_1\rt| dH^1 z\nn\\
&\overset{(\ref{peq25})}{\leq}& 20000 h_0 \times (\cfc)^{\frac{Q}{2p}}
\lt(\frac{Q}{\mu}\rt)^{\frac{22Q}{p}} C_p^{\frac{1}{p}}\ep^{\frac{p}{960 Q^2}}\nn\\
&\leq& h_0 10^{\frac{36Q}{p}}\lt(\frac{Q}{\mu}\rt)^{\frac{22Q}{p}} C_p^{\frac{1}{p}}\ep^{\frac{p}{960 Q^2}}.
\end{eqnarray}
Let 
\begin{equation}
\label{breq220}
\varpi=10^{\frac{36Q}{p}}
\lt(\frac{Q}{\mu}\rt)^{\frac{22Q}{p}} C_p^{\frac{1}{p}}\text{ and }\beta=\frac{p}{960 Q^2}.
\end{equation}
Note 
\begin{equation}
\label{oleq1}
h_0\overset{(\ref{breq80}) }{\geq} \frac{\mu^{\frac{21 Q}{p}} Q^{-\frac{20 Q}{p}} C_p^{-\frac{1}{2p}}}{96 \pi \lt(9(1352)^{18 Q}\rt)^{\frac{1}{2p}}}. 
\end{equation}
Now let 
\begin{equation}
\label{peq27}
w(z)=\phi_u'(z)-\zeta_1 \phi_v'(z). 
\end{equation}
Hence by Cauchy's integral formula we have that 
\begin{eqnarray}
\label{peq191}
\lt|w^{(k)}(x_0)\rt|&=&\frac{k!}{2\pi}\int_{\partial B_q(x_0)} \lt|\frac{w(\zeta)}{(\zeta-x_0)^{k+1}}\rt| d\zeta\nn\\
&\leq & \frac{k!}{2\pi q^{k+1}}\int_{\partial B_q(x_0)} \lt|w(\zeta)\rt| d\zeta\nn\\
&\overset{(\ref{breq220}),(\ref{jejq1})}{\leq}& \frac{k!}{2\pi q^{k+1}} \varpi h_0 \ep^{\beta}\nn\\
&\overset{(\ref{oopeq70})}{\leq}& \frac{2 k! \varpi \ep^{\beta}}{q^k}.
\end{eqnarray}
By the local Talyor Theorem we have 
\begin{equation}
\label{eq23}
w(z)=\sum_{k=0}^{m} \frac{w^{(k)}(x_0)}{k!}(z-x_0)^k+(z-x_0)^{m+1}w_m(z)
\end{equation}
where $w_m(z)=\frac{1}{2\pi i}\int_{\partial B_{\frac{3\mu}{2}}(x_0)} \frac{w(\zeta)}{(\zeta-x_0)^m(\zeta-z)} d\zeta$ for 
any $z\in B_{\frac{3\mu}{2}}(x_0)$. Hence for $z\in B_{\mu}(x_0)$
\begin{eqnarray}
\label{oppeq11}
\lt|w_m(z)\rt|&\leq& \frac{1}{2\pi}\int_{\partial B_{\frac{3\mu}{2}}(x_0)} \frac{\lt|w(\zeta)\rt|}{\lt|\zeta-x_0\rt|^m\lt|\zeta-z\rt|} dz\nn\\
&\overset{(\ref{peq27}),(\ref{eq3})}{\leq}& \frac{16}{\mu} \int_{\partial B_{\frac{3\mu}{2}}(x_0)} \frac{1}{(\frac{3\mu}{2})^m \frac{\mu}{2}}\nn\\
&\leq& 64\pi \mu^{-2}\lt(\frac{3\mu}{2}\rt)^{1-m}.
\end{eqnarray}
So for any $z\in B_{\mu}(x_0)$ we have 
\begin{eqnarray}
\label{eqb1}
\lt|w(z)\rt|&\overset{(\ref{eq23}),(\ref{oppeq11})}{\leq}& \sum_{k=0}^m \frac{\lt|w^{(k)}(x_0)\rt|}{k!}\lt|z-x_0\rt|^k
+\lt|z-x_0\rt|^{m+1} 64\pi \mu^{-2} \lt(\frac{3\mu}{2}\rt)^{1-m}\nn\\
&\overset{(\ref{peq191})}{\leq}& 2\sum_{k=0}^m \varpi \ep^{\beta} \lt(\frac{\mu}{q}\rt)^k
+64\pi\lt(\frac{3}{2}\rt)^{1-m}.
\end{eqnarray}
Let 
\begin{eqnarray}
\label{breq81}
\alpha&=&\frac{h_0}{\mu}\overset{(\ref{breq80})}{=}\frac{\mu}{96\pi}\lt(\frac{\mu^{38Q}}{9 C_p (1352)^{18 Q} Q^{40Q}}\rt)^{\frac{1}{2p}}\nn\\
&\geq& \frac{\mu^{\frac{20Q}{p}}}{(96\pi\times 3 \times (1352)^9)^{\frac{Q}{p}}C_p Q^{20Q}}\geq
\frac{\mu^{\frac{20Q}{p}}}{C_p Q^{20Q} 10^{\frac{32 Q}{p}}}.
\end{eqnarray}
Now note that $q\in (\frac{\alpha \mu}{2},\alpha \mu)$, 
$\lt(\frac{\mu}{q}\rt)^k\leq \lt(\frac{2}{\alpha}\rt)^k$. So note as 
\begin{equation}
\label{peqq1}
\alpha<1.
\end{equation}
Hence 
\begin{equation}
\label{peq40}
\sum_{k=0}^m \lt(\frac{2}{\alpha}\rt)^k\leq  \frac{\lt(\frac{2}{\alpha}\rt)^{m+1}}{\frac{2}{\alpha}-1}\leq \lt(\frac{2}{\alpha}\rt)^{m+1}.
\end{equation}
Thus
\begin{eqnarray}
\label{ecz2}
\sum_{k=0}^m  \varpi \ep^{\beta} \lt(\frac{\mu}{q}\rt)^k&\leq&   \varpi \ep^{\beta} 
\sum_{k=0}^m \lt(\frac{2}{\alpha}\rt)^k\nn\\
&\overset{(\ref{peq40})}{\leq}& \varpi \ep^{\beta} \lt(\frac{2}{\alpha}\rt)^{m+1}.
\end{eqnarray}
So
\begin{equation}
\label{eqd1}
\lt|w(z)\rt|\overset{(\ref{eqb1}),(\ref{ecz2})}{\leq} 
2\varpi \ep^{\beta} \lt(\frac{2}{\alpha}\rt)^{m+1}+64\pi\lt(\frac{3}{2}\rt)^{1-m}.
\end{equation}
Let $m$ be the smallest integer such that 
\begin{equation}
\label{eqa1}
2\varpi\ep^{\beta}\lt(\frac{2}{\alpha}\rt)^{m}\geq 
64\pi\lt(\frac{2}{3}\rt)^m.
\end{equation}
So
\begin{eqnarray}
\label{oppeq80}
\ep^{\beta}&\geq& \frac{32\pi}{\varpi}\lt(\frac{\alpha}{2}\rt)^{m}\lt(\frac{2}{3}\rt)^m\nn\\
&=& 
\lt(\lt( \frac{32\pi}{\varpi}\rt)^{\frac{1}{m}} \lt(\frac{\alpha}{3}\rt)\rt)^m.
\end{eqnarray}
Thus as $\varpi\overset{(\ref{breq220})}{>}32 \pi$ and $\alpha\overset{(\ref{peqq1})}{<}1$
\begin{eqnarray}
\lt|\log(\ep^{\beta})\rt|&\overset{(\ref{oppeq80})}{\leq}& \lt|\log\lt(\lt(\lt(\frac{32\pi}{\varpi}\rt)^{\frac{1}{m}}
\lt(\frac{\alpha}{3}\rt)\rt)^m\rt) \rt|\nn\\&=&m\lt|\log\lt(\lt(\frac{32\pi}{\varpi}\rt)^{\frac{1}{m}} \lt(\frac{\alpha}{3}\rt)\rt) \rt|.\nn
\end{eqnarray}
Hence 
\begin{equation}
\label{eqa2}
\frac{\log\lt(\ep^{\beta}\rt)}{\log\lt(\lt(\frac{32\pi}{\varpi}\rt)^{\frac{1}{m}}  
\lt(\frac{\alpha}{3}\rt)\rt)}\leq m.
\end{equation}
So 
\begin{eqnarray}
\label{eqa5}
96\pi\lt(\frac{2}{3}\rt)^m&\leq& 96\pi \lt(\frac{2}{3}\rt)^{\frac{\log\lt(\ep^{\beta}\rt)}
{\log\lt(\lt(\frac{32\pi}{\varpi}\rt)^{\frac{1}{m}} \lt(\frac{\alpha}{3}\rt)\rt)}}\nn\\
&=&96\pi\lt(e^{\log\lt(\frac{2}{3}\rt)}\rt)^{\frac{\log\lt(\ep^{\beta}\rt)}{\log\lt(\lt(\frac{32\pi}{\varpi}\rt)^{\frac{1}{m}}  \lt(\frac{\alpha}{3}\rt)\rt)}}\nn\\
&=&96\pi \lt(\ep^{\beta}\rt)^{\frac{\log\lt(\frac{2}{3}\rt)}{\log\lt(\lt(\frac{32\pi}{\varpi}\rt)^{\frac{1}{m}} \lt(\frac{\alpha}{3}\rt)\rt)}}.
\end{eqnarray}
Using the fact that $\varpi\overset{(\ref{breq220})}{>}32 \pi$ so $1>\frac{32 \pi}{\varpi}$ and thus 
$1>\lt(\frac{32 \pi}{\varpi}\rt)^{\frac{1}{m}}>\lt(\frac{32 \pi}{\varpi}\rt)$. Thus we have 
$$
1\overset{(\ref{peqq1})}{\geq} \lt(\frac{32\pi}{\varpi}\rt)^{\frac{1}{m}}\lt(\frac{\alpha}{3}\rt) 
\geq \frac{10\pi\alpha}{\varpi}, 
$$
so 
\begin{equation}
\label{oopeq35}
\lt|\log\lt(\lt(\frac{32\pi}{\varpi}\rt)^{\frac{1}{m}}\lt(\frac{\alpha}{3}\rt)\rt)\rt|\leq \lt|\log\lt(\frac{10\pi\alpha}{\varpi}\rt)\rt|.
\end{equation}
So we have 
\begin{equation}
\label{eqa6}
96\pi\lt(\frac{2}{3}\rt)^m\overset{(\ref{eqa5})}{\leq} 96\pi\ep^{\frac{\beta\log(\frac{2}{3})}
{\log\lt(\lt(\frac{32\pi \alpha}{\varpi}\rt)^{\frac{1}{m}}\lt(\frac{\alpha}{3}\rt)\rt) }}
\overset{(\ref{oopeq35})}{\leq}  
96\pi\ep^{\frac{-\beta}{3\log\lt(\frac{10\pi \alpha}{\varpi}\rt)}}  .
\end{equation}
Since $m$ is the smallest integer such that (\ref{eqa1}) holds true we have 
\begin{equation}
\label{ecz10}
2\varpi\ep^{\beta}\lt(\frac{2}{\alpha}\rt)^{m-1} \overset{(\ref{eqa1})}{\leq} 64\pi \lt(\frac{2}{3}\rt)^{m-1}= 96 \pi \lt(\frac{2}{3}\rt)^m\overset{(\ref{eqa6})}{\leq} 
96\pi \ep^{\frac{-\beta}{3\log\lt(\frac{10\pi \alpha}{\varpi}\rt)}}.
\end{equation}
Thus as $\alpha^2\overset{(\ref{breq81})}{\geq} \frac{\mu^{\frac{40 Q}{p}}}{C_p^2 Q^{40Q} 10^{\frac{64Q}{p}}}$ 
\begin{eqnarray}
\label{nika1}
\varpi\ep^{\beta}\lt(\frac{2}{\alpha}\rt)^{m+1} &=&\frac{4}{\alpha^2}\varpi\ep^{\beta}\lt(\frac{2}{\alpha}\rt)^{m-1}\nn\\
&\overset{(\ref{ecz10})}{\leq}&\frac{2}{\alpha^2}
\times 96 \pi \ep^{\frac{-\beta}{3\log\lt(\frac{10\pi \alpha}{\varpi}\rt)}}\nn\\
&\leq& \frac{2\times 96 \pi \times 10^{\frac{64Q}{p}}  C_p^2 Q^{40Q}}{\mu^{\frac{40Q}{p}}}
\ep^{\frac{-\beta}{3\log\lt(\frac{10\pi \alpha}{\varpi}\rt)}}\nn\\
&\leq&C_p^2 \CI_0 \ep^{\frac{-\beta}{3\log\lt(\frac{10\pi \alpha}{\varpi}\rt)}},
\end{eqnarray}
where $\CI_0=\CI_0(p,Q)$.

So putting (\ref{nika1}) and (\ref{eqa6}) together with (\ref{eqd1}) we have 
\begin{eqnarray}
\lt|w(z)\rt|&\leq& C_p^2 \CI_0 \ep^{\frac{-\beta}{3\log\lt(\frac{10\pi \alpha }{\varpi}\rt)}}
+96\pi \ep^{\frac{-\beta}{3\log\lt(\frac{10\pi \alpha}{\varpi}\rt)}} .    \nn\\
&\leq&   C_p^2 \CI_1  \ep^{\frac{-\beta}{3\log\lt(\frac{10\pi \alpha}{\varpi}\rt)}}   \text{ for all }z\in B_{\mu}(z_0).\nn
\end{eqnarray}
Hence 
\begin{equation}
\label{oopeq91}
\|\phi'_u-\zeta_1\phi_v'\|_{L^{\infty}(B_{\mu}(z_0))}\leq   
C_p^2 \CI_1 \ep^{\frac{-\beta}{3\log\lt(\frac{10\pi \alpha}{\varpi}\rt)}}
=C_p^2 \CI_1 \ep^{\frac{\beta}{3\log\lt(\frac{\varpi}{10\pi \alpha}\rt)}}. 
\end{equation}
Now 
\begin{eqnarray}
\frac{10\pi\alpha}{\varpi}&\overset{(\ref{breq81})}{\geq}&\frac{10\pi}{\varpi} \frac{\mu^{\frac{20Q}{p}}}{C_p Q^{20} 10^{\frac{32 Q}{p}}}\nn\\
&\geq&\lt(\frac{\mu^2}{2000 \varpi C_p Q^2}\rt)^{\frac{10 Q}{p}}.\nn
\end{eqnarray}
So 
\begin{eqnarray}
\label{lleq55}
\frac{\varpi}{10 \pi \alpha}&\leq& \lt(\frac{2000 \varpi C_p Q^2}{\mu^2}\rt)^{\frac{10 Q}{p}}\nn\\
&\leq& (2000)^{\frac{10Q}{p}} \varpi^{\frac{10 Q}{p}}C_p^{\frac{10Q}{p}}\lt(\frac{Q}{\mu}\rt)^{\frac{20 Q}{p}}\nn\\
&\overset{(\ref{breq220}),(\ref{oleq1})}{\leq}& C_p^{\frac{20 Q}{p}} (2000)^{\frac{10Q}{p}} 10^{\frac{360Q}{p}}\lt(\frac{Q}{\mu}\rt)^{\frac{240 Q^2}{p^2}}.
\end{eqnarray}
Now 
\begin{equation}
\label{ooleq42}
\mu\overset{(\ref{peq51})}{\geq} (10^{10}\times 2^7 Q^6)^{-12Q}. 
\end{equation}
Thus 
\begin{eqnarray}
\label{oppeq30}
\frac{\varpi}{10 \pi \alpha}&\overset{(\ref{lleq55})}{\leq}& C_p^{\frac{20 Q}{p}}\lt(2000  \times 10^{36} \rt)^{\frac{10Q^2}{p^2}}
\times Q^{\frac{240 Q^2}{p^2}} (10^{10}\times 2^7 Q^6)^{\frac{2880 Q^3}{p^2}}\nn\\
&\leq& C_p^{\frac{20 Q}{p}}\lt(2000\times 10^{36} \rt)^{\frac{10Q^2}{p^2}}\times \lt(10^{10}\times 2^7 \rt)^{\frac{2880 Q^3}{p^2}}Q^{\frac{240 Q^2}{p^2}+\frac{17280 Q^3}{p^2}}\nn\\
&\leq& C_p^{\frac{20 Q}{p}}\lt((2000 \times 10^{36})^{10}\times (10^{10}\times 2^7)^{2880} \rt)^{\frac{Q^3}{p^2}}
Q^{17520\frac{Q^3}{p^2}}\nn\\
&\leq& (10 C_p Q)^{35262\frac{Q^3}{p^2}}
\end{eqnarray}
So 
\begin{eqnarray}
\label{oppeq31}
\ep^{\frac{\beta}{3 \log\lt(\frac{\varpi}{10\pi \alpha}\rt)}}&\overset{(\ref{oppeq30})}{\leq}&
\ep^{\frac{p^2 \beta}{35262 Q^3 \log(10 C_p Q)}}\nn\\
&\overset{(\ref{breq220})}{\leq}& \ep^{\frac{p^3}{35262  \times 960 Q^5 \log(10 C_p Q)}}\nn\\
&\leq& \ep^{\frac{p^3}{4\times 10^7 Q^5 \log(10 C_p Q)}}
\end{eqnarray}
Thus as $x_0\in B_{\frac{\mu}{2}}(w_u(0))$ we know 
$B_{\frac{\mu}{2}}(w_u(0))\subset B_{\mu}(x_0)$
\begin{equation}
\|\phi_u'-\zeta_1 \phi'_v\|_{L^{\infty}(B_{\mu}(x_0))}\overset{(\ref{oopeq91}),(\ref{oppeq30}),(\ref{oppeq31})}{\leq} C_p \CI_1 \ep^{\frac{p^3}{4\times 10^{7} Q^5 \log(10 C_p Q)}}
\end{equation}
and hence we have established (\ref{eq60}).\qd\qd$\Box$\nl\nl
%
%
%
%
%
\subsection{Proof of Proposition \ref{P1} completed}
Now 
$$
D u(x)=D \phi_u(w_u(x))D w_u(x)\text{ and }D v(x)=D \phi_v(w_v(x))D w_v(x).
$$ 
So
\begin{eqnarray}
\label{eqw400}
\int_{B_{\gamma}} \lt|Du(x)-R Dv(x)\rt| dx&=&\int_{B_{\gamma}} \lt|D\phi_u(w_u(x))D w_u(x)-R D\phi_v(w_v(x))D w_v(x) \rt| dx\nn\\
&\leq& 
\int_{B_{\gamma}} \lt|\lt(D\phi_u(w_u(x))-R D\phi_v(w_u(x))\rt) D w_u(x) \rt| dx\nn\\
&~&+\int_{B_{\gamma}}\lt|D\phi_v(w_u(x))(Dw_u(x)- Dw_v(x))\rt| dx\nn\\
&~&+\int_{B_{\gamma}}\lt|(D\phi_v(w_u(x))- D\phi_v(w_v(x)))D w_v(x)\rt| dx.\nn
\end{eqnarray}
So to deal with the last term 
\begin{eqnarray}
\label{peq55}
&~&\int_{B_{\gamma}} \lt|(D\phi_v(w_u(x))- D\phi_v(w_v(x)))D w_v(x)\rt| dx\nn\\
&~&\qd\qd\qd\overset{(\ref{eq5}),(\ref{ezq1}),(\ref{pw2})}{\leq}
\frac{192\pi}{\mu^2}\int_{B_{\gamma}} \lt|w_u(x)-w_v(x)\rt|\lt|D w_v(x)\rt| dx\nn\\
&~&\qd\qd\qd\overset{(\ref{eqz40})}{\leq} \frac{192 \pi}{\mu^2} \times 1104 Q^2 \pi C_p \ep^{\frac{p}{120 Q^2}}\sqrt{\pi}
\lt(\int_{B_{\gamma}} \lt|D w_u(x)\rt|^2 dx \rt)^{\frac{1}{2}}\nn\\
&~&\qd\qd\qd\overset{(\ref{eaz1})}{\leq} \frac{\CI_2 Q^3 C_p}{\mu^2} \ep^{\frac{p}{120 Q^2}}.\nn
\end{eqnarray}
And
\begin{eqnarray}
\int_{B_{\gamma}} \lt|D\phi_v(w_u(x))\lt(D w_u(x)-D w_v(x)\rt)\rt| dx &\overset{(\ref{eq3}),(\ref{equa127})}{\leq}& \frac{32}{\mu}\times 24 \pi^2 C_p \sqrt{Q}\ep^{\frac{p}{48 Q^2}}\nn\\
&\leq& \CI_2 C_p \sqrt{Q}\ep^{\frac{p}{48 Q^2}}.
\end{eqnarray}
So 
\begin{eqnarray}
\int_{B_{\gamma}} \lt|Du-R Dv\rt| dx&=&\sqrt{\pi}\lt(\int_{B_{\gamma}} \lt|D \phi_u(w_u(x))-R D\phi_v(w_u(x))\rt|^2\lt|D w_u(x)\rt|^2 dx\rt)^{\frac{1}{2}}\nn\\
&~&\qd\qd+ \frac{\CI_2 Q^3 C_p}{\mu^2} \ep^{\frac{p}{120 Q^2}}+\CI_2 C_p \sqrt{Q} \ep^{\frac{p}{48 Q^2}}\nn\\
&\leq& \sqrt{Q\pi}\lt(\int_{B_{\gamma}} \lt|D \phi_u(w_u(x))-R D\phi_v(w_u(x))\rt|^2\det(D w_u(x))\; dx\rt)^{\frac{1}{2}}\nn\\ 
&~&\qd\qd+ \frac{\CI_2 Q^3 C_p}{\mu^2} \ep^{\frac{p}{120 Q^2}}\nn\\
&\overset{(\ref{ezq1})}{\leq}& \sqrt{Q\pi}\lt(\int_{B_{\frac{\mu}{2}}(w_u(0))} \lt|D\phi_u(z)-R D\phi_v(z)\rt|^2 dz \rt)^{\frac{1}{2}}
+\frac{\CI_2 Q^3 C_p}{\mu^2} \ep^{\frac{p}{120 Q^2}}\nn\\
&\overset{(\ref{eq60})}{\leq}& C_p \CI_4 \ep^{\frac{p^3}{4\times 10^7 Q^5 \log(10 C_p Q)}}.\qd\qd\Box
\end{eqnarray}

%
%
%

\section{Proof of Theorem \ref{T1}}

Let $\ti{u}(z)=\frac{u(z)}{4}$ and $\ti{v}(z)=\frac{v(z)}{4}$. So 
\begin{equation}
\label{lleq11}
\int_{B_1} \lt|D\ti{u}\rt| dz\leq \frac{1}{4}
\end{equation}
and 
\begin{eqnarray}
\label{fgfeq1}
\int_{B_1} \det(D\ti{u})^{-p} dx&=&16^p\int_{B_1} \det(D u)^{-p} dx \nn\\
&=& 16^p C_p.
\end{eqnarray}
Note also from (\ref{fg1})
\begin{equation}
\label{lleq5}
\int_{B_1} \lt|S(D\ti{u}(z))-S(D \ti{v}(z))\rt|^2 dz\leq \ep.
\end{equation}

\em Step 1. \rm For any set $S\subset B_1$ with $\lt|S\rt|>0$ we will show 
\begin{equation}
\label{peq701}
\int_{S} \det(D\ti{u}(z)) dz\geq 16^{-1} C_p^{-\frac{1}{p}}\lt|S\rt|^{\frac{2-p}{p-p^2}}.
\end{equation}

\em Proof of Step 1. \rm Note 
\begin{eqnarray}
\label{oopeq96}
\lt|S\rt|&=&\int_{S} \det(D\ti{u}(z))^{\frac{p}{2}}\det(D\ti{u}(z))^{-\frac{p}{2}} dz\nn\\
&\overset{(\ref{fgfeq1})}{\leq}&\lt( \int_{S} \det(D\ti{u}(z))^{p} dz\rt)^{\frac{1}{2}}16^{\frac{p}{2}}\sqrt{C_p}.
\end{eqnarray}

Let $q=\frac{1}{p}$, $q'=\frac{q}{q-1}=\frac{p^{-1}}{p^{-1}-1}=\frac{p}{1-p}$. So 
using Holder's inequality for the second inequality  
\begin{eqnarray}
\label{oopeq97}
\lt(\frac{\lt|S\rt|}{16^{\frac{p}{2}}\sqrt{C_p}}\rt)^2&\overset{(\ref{oopeq96})}{\leq}& 
\int_{S} \det(Du(z))^p dz\nn\\
&=&\lt(\int_{S} \det(Du(z))^{pq} dz\rt)^{\frac{1}{q}}\lt|S\rt|^{\frac{1}{q'}}\nn\\
&=&\lt(\int_{S} \det(Du(z)) dz\rt)^{p}\lt|S\rt|^{\frac{1-p}{p}}.
\end{eqnarray}
So 
\begin{eqnarray}
16^{-p}C_p^{-1}\lt|S\rt|^{\frac{2-p}{1-p}}&=&16^{-p}C_p^{-1} \lt|S\rt|^2 \lt|S\rt|^{\frac{p}{1-p}}\nn\\
&\overset{(\ref{oopeq97})}{\leq}& \lt( \int_{S} \det(Du(z)) dz \rt)^p.
\end{eqnarray}
Thus 
$$
16^{-1} C_p^{-\frac{1}{p}}\lt|S\rt|^{\frac{2-p}{p-p^2}}\leq \int_{S} \det(Du(z)) dz
$$
so we have established (\ref{peq701}).\nl

\em Step 2. \rm Let $\lt\{B_{\frac{\gamma}{2}}(x_k):k=1,2,\dots N\rt\}$ be collection 
such that 
\begin{equation}
\label{peq802}
\sum_{k=1}^N \chara_{B_{\frac{\gamma}{2}}(x_k)}\leq 5
\end{equation}
and 
\begin{equation}
\label{peq803}
B_{\frac{1}{2}}\subset \bigcup_{k=1}^N B_{\frac{\gamma}{2}}(x_k).
\end{equation}
We will order these balls so that $B_{\frac{\gamma}{2}}(x_k)\cap B_{\frac{\gamma}{2}}(x_{k+1})\not=\emptyset$ 
for $k=1,2,\dots N-1$. 

Let $u_k(z)=2\ti{u}(x_k+\frac{z}{2})$ and $v_k(z)=2\ti{v}(x_k+\frac{z}{2})$. Note 
$$
\int_{B_1} \lt|D u_k(z)\rt| dz=\int_{B_1} \lt|D \ti{u}\lt(x_k+\frac{z}{2}\rt)\rt| dz
\leq 4 \int_{B_{\frac{1}{2}}(x_k)} \lt|D \tilde{u}(z)\rt| dz\overset{(\ref{lleq11})}{\leq} 1. 
$$
And 
\begin{equation}
\label{lleq6}
\int_{B_1} \det\lt(D \tilde{u}\lt(x_k+\frac{z}{2}\rt)\rt)^{-p} dz
\leq 4\int_{B_1} \det\lt(D \tilde{u}(y)\rt)^{-p} dy\overset{(\ref{fgfeq1})}{\leq} 64 C_p.
\end{equation}
Note also 
\begin{eqnarray}
\int_{B_1} \lt|S(Du_k)-S(Dv_k)\rt|^2 dz&\leq&  \int_{B_1} \lt|S\lt(D\ti{u}\lt(x_k+\frac{x}{2}\rt)\rt)
-S\lt(D\ti{u}\lt(x_k+\frac{x}{2}\rt)\rt) \rt|^2 dx\nn\\
&=&4 \int_{B_1} \lt|S\lt(D\ti{u}(y)\rt)-S\lt(D\ti{v}(y)\rt) \rt| dy\nn\\
&\overset{(\ref{lleq5})}{\leq}& 4\ep.  
\end{eqnarray}
So we can apply Proposition \ref{P1} and for some $R_k\in SO(2)$ we have 
$$
\int_{B_{\gamma}} \lt|Dv_k(z)-R_k Du_k(z)\rt| dz\leq \CI_4 \CI_p \ep^{\frac{p^3}{4\times 10^7 Q^5 \log(10 C_p Q)}}
$$
We will show that 
\begin{equation}
\label{lleq8}
\lt|R_1-R_k\rt|\leq c \gamma^{-1-\frac{(2-p)}{p-p^2}} C_p^{\frac{1}{p}} \ep^{\frac{p^3}{10^8 Q^5 \log(10 C_p Q)}}\text{ for }k=1,2,\dots N-1.
\end{equation}

\em Proof of Step 2. \rm The existence of a collection 
$\lt\{B_{\frac{\gamma}{2}}(x_1), B_{\frac{\gamma}{2}}(x_2), \dots   B_{\frac{\gamma}{2}}(x_N) \rt\}$ satisfying 
(\ref{peq802}), (\ref{peq803}) follows by the $5r$ covering theorem.

Rescaling $v_k$ and $u_k$ we have 
\begin{equation}
\label{lleq1}
\int_{B_{\gamma}(x_k)} \lt| D\tilde{v}(z)-R_k D \tilde{u}(z)\rt| dz\leq \CI_4 C_p \ep^{\frac{p^3}{4\times 10^7 Q^5\log(10 C_p Q)}}.
\end{equation}
So 
\begin{equation}
\label{lleq2}
\int_{B_{\gamma}(x_k)\cap B_{\gamma}(x_{k+1})} \lt|(R_k-R_{k+1}) D\tilde{u}(z)\rt| dz\leq \CI_4 C_p \ep^{\frac{p^3}{4\times 10^7 Q^5 \log(10 C_p Q)}}.
\end{equation}
Let 
\begin{equation}
\label{eqll2}
\BI_1=\lt\{z: \lt|D\tilde{u}(z)\rt|>2\gamma^{-2}\rt\}.
\end{equation}
So $\lt|\BI_1\rt|\overset{(\ref{lleq11})}{\leq} \frac{\gamma^2}{8}$. Since $B_{\frac{\gamma}{2}}(x_k)\cap B_{\frac{\gamma}{2}}(x_{k+1})\not = \emptyset$ for 
$k=1,2,\dots N-1$. So 
$$
\lt|B_{\gamma}(x_k)\cap B_{\gamma}(x_{k+1})\rt|\geq \frac{\gamma^2}{4}\text{ for }k=1,2,\dots N-1.  
$$
Thus  
\begin{equation}
\label{lleq2.5}
\lt| B_{\gamma}(x_k)\cap B_{\gamma}(x_{k+1})\backslash \BI_1\rt|\geq \frac{\gamma^2}{8}\text{ for }k=1,2,\dots N-1. 
\end{equation}
Now 
\begin{eqnarray}
\label{lleq112}
&~&\int_{B_{\gamma}(x_k)\cap B_{\gamma}(x_{k+1})\backslash \BI_1} 
\det(R_k-R_{k+1})\det(D \ti{u}(z)) dz\nn\\
&~&\qd\qd\qd= \int_{B_{\gamma}(x_k)\cap B_{\gamma}(x_{k+1})\backslash \BI_1} 
\det\lt((R_k-R_{k+1})D \ti{u}(z)\rt) dz\nn\\
&~&\qd\qd\qd\leq \int_{B_{\gamma}(x_k)\cap B_{\gamma}(x_{k+1})\backslash \BI_1} 
\|(R_k-R_{k+1})D \ti{u}(z)\|^2 dz\nn\\
 &~&\qd\qd\qd\leq 2 \int_{B_{\gamma}(x_k)\cap B_{\gamma}(x_{k+1})\backslash \BI_1} 
\|D \ti{u}(z)\|\lt|(R_k-R_{k+1})D \ti{u}(z)\rt| dz\nn\\
 &~&\qd\qd\qd\overset{(\ref{eqll2})}{\leq}  4 \gamma^{-2}\int_{B_{\gamma}(x_k)\cap B_{\gamma}(x_{k+1})\backslash \BI_1}  
\lt|(R_k-R_{k+1})D \ti{u}(z)\rt| dz\nn\\
 &~&\qd\qd\qd\overset{(\ref{lleq2})}{\leq}   c \gamma^{-2} C_p \ep^{\frac{p^3}{4\times 10^7 Q^5 \log(10 C_p Q)}}.
\end{eqnarray}
Hence 
\begin{eqnarray}
\det(R_k-R_{k+1})C_p^{-\frac{1}{p}}\gamma^{\frac{2(2-p)}{p-p^2}}&\overset{(\ref{lleq2.5})}{\leq}& c\det(R_k-R_{k+1})C_p^{-\frac{1}{p}}
\lt|B_{\gamma}(x_k)\cap B_{\gamma}(x_{k+1})\backslash \BI_1 \rt|^{\frac{2-p}{p-p^2}}\nn\\
&\overset{(\ref{peq701})}{\leq}& \int_{B_{\gamma}(x_k)\cap B_{\gamma}(x_{k+1})\backslash \BI_1 } \det\lt(R_k-R_{k+1}\rt) \det(D \ti{u}(z)) dz\nn\\
&\overset{(\ref{lleq112})}{\leq}& c \gamma^{-2} C_p \ep^{\frac{p^3}{4\times 10^7 Q^5 \log(10 C_p Q)}}.\nn
\end{eqnarray}
Thus 
\begin{equation}
\label{lleq6.7}
\det(R_k-R_{k+1})\leq c \gamma^{-2-\frac{2(2-p)}{p-p^2}} C_p^{\frac{2}{p}} \ep^{\frac{p^3}{4\times 10^7 Q^5 \log(10 C_p Q)}}.
\end{equation}

Note that if $R_{\alpha}=\lt(\begin{smallmatrix}\cos(\alpha) & -\sin(\alpha) \\ \sin(\alpha) & \cos(\alpha)\end{smallmatrix}\rt)$, $R_{\beta}=\lt(\begin{smallmatrix}\cos(\beta) & -\sin(\beta) \\ \sin(\beta) & \cos(\beta)\end{smallmatrix}\rt)$ then $\det(R_{\alpha}-R_{\beta})=2(1-\cos(\alpha-\beta))$ thus from (\ref{lleq6.7}) we have 
 \begin{equation}
\label{lleq7}
\lt|R_k-R_{k+1}\rt|\leq c \gamma^{-1-\frac{(2-p)}{p-p^2}} C_p^{\frac{1}{p}} \ep^{\frac{p^3}{10^8 Q^5 \log(10 C_p Q)}}\text{ for }k=1,2,\dots N-1.
\end{equation}
So we have established (\ref{lleq8}). \qd\qd $\Box$\nl

\subsection{Proof of Theorem \ref{T1} completed} 

\begin{eqnarray}
\int_{B_{\gamma}(x_k)} \lt| D \ti{v}(z)-R_1 D \ti{u}(z)\rt| dz&\leq& 
\int_{B_{\gamma}(x_k)} \lt| D \ti{v}(z)-R_k D \ti{u}(z)\rt| dz+\lt|(R_k-R_1) D\ti{u}(z)\rt| dz\nn\\
&\overset{(\ref{lleq1}),(\ref{lleq8})}{\leq}& \CI_4 C_p \ep^{\frac{p^3}{10^8 Q^5\log(10 C_p Q)}}
+ \CI_5 C_p^{\frac{1}{p}} \ep^{\frac{p^3}{ 10^8 Q^5 \log(10 C_p Q)}}\int_{B_{\gamma}(x_k)} \lt| D \ti{u}\rt| dz. \nn
\end{eqnarray}
Thus 
\begin{eqnarray}
\int_{B_{\frac{1}{2}}} \lt| D \ti{v}- R_1 D \ti{u}\rt| dz&\overset{(\ref{peq802}),(\ref{peq803})}{\leq}&c\sum_{k=1}^N  C_p^{\frac{1}{p}}\ep^{\frac{p^3}{10^8 Q^5 \log(10 C_p Q)}}\int_{B_{\gamma}(x_k)} \lt| D \ti{u}\rt| dz
+ c C_p \ep^{\frac{p^3}{10^8 Q^5 \log(10 C_p Q)}}\nn\\
&\overset{(\ref{peq802}),(\ref{lleq11})}{\leq}&  c C_p^{\frac{1}{p}} \ep^{\frac{p^3}{10^8 Q^5 \log(10 C_p Q)}}.
\end{eqnarray}
Rescaling gives (\ref{fest1}) and this completes the proof of Theorem \ref{T1}.\qd\qd $\Box$

\section{Examples}
\label{example}

We can show that any estimate has to lose at least a root power. \nl

\em Example 1. \rm

Let $f(z)=\frac{z^{k+1}}{k+1}$, $g(z)=\frac{z^{k+2}}{k+2}$. So rewriting these functions as vector valued functions 
of two variables we have 
\begin{equation}
\label{lleq110}
D\ti{f}(x,y)=\lt[z^k\rt]_M\text{ and }D\ti{g}(x,y)=\lt[z^{k+1}\rt]_M.
\end{equation}
Now 
\begin{equation}
\lt[z^k\rt]_M=\lt|z\rt|^k  \lt(\begin{array}{cc} \cos(k \mathrm{Arg}(z)) &  -\sin(k \mathrm{Arg}(z)) 
\\  \sin(k \mathrm{Arg}(z)) & \cos(k \mathrm{Arg}(z))   \end{array}\rt)
\end{equation}
and 
\begin{equation}
\lt[z^{k+1}\rt]_M=\lt|z\rt|^{k+1}  \lt(\begin{array}{cc} \cos((k+1) \mathrm{Arg}(z)) &  -\sin((k+1) \mathrm{Arg}(z)) 
\\  \sin((k+1) \mathrm{Arg}(z)) & \cos((k+1) \mathrm{Arg}(z))   \end{array}\rt)
\end{equation}
Thus 
$$
\mathrm{Sym}(D\ti{f}(x,y))=(x^2+y^2)^{\frac{k}{2}}Id\text{ and }\mathrm{Sym}(D\ti{g}(x,y))=(x^2+y^2)^{\frac{k+2}{2}}Id.
$$
So note 
\begin{eqnarray}
\int_{B_1} \lt|\mathrm{Sym}(D\ti{f})- \mathrm{Sym}(D\ti{g})\rt| dx&=& 
\int_{0}^1 \int_{\partial B_r} \lt|r^{k}-r^{k+1}\rt| dH^1 z dr \nn\\
&=& 2\pi \lt(\frac{1}{k+1}-\frac{1}{k+2}\rt)\nn\\
&=& \frac{2\pi}{(k+1)(k+2)}.\nn
\end{eqnarray}
A slightly longer calculation shows that 
\begin{equation}
\label{eqddg1}
\int_{B_1} \lt|D\ti{f}-R_{\theta}D\ti{g}\rt| dz\geq \frac{c}{k}\text{ for any }\theta\in \lt(0,2\pi\rt].
\end{equation}

%
%
%

\begin{a3}
\label{CJE1}
There exists a sequence of positive numbers $\ep_k\rightarrow 0$ and a sequence of pairs 
of $Q$-Quasiregular maps $u_k:B_1\rightarrow \R^2$, $v_k:B_1\rightarrow \R^2$ with 
$\int_{B_1} \lt| Du_k\rt|^2 dz\leq 1$ such that 
$$
\int_{B_1} \lt|S(Du_k)-S(Dv_k)\rt|^2 dz=\ep_k
$$
and 
$$
\int_{B_{\frac{1}{2}}} \lt|Du_k-R_{\theta} Dv_k\rt| dz\geq 1\text{ for all }R_{\theta}\in SO(2).
$$
\end{a3}
\em Sketch of proof of Conjecture \ref{CJE1}. \rm Let $k$ be a large integer. Let $w_m=e^{\frac{2\pi i m}{k}}$. 
A natural approach is to define function 
\begin{equation}
h(z):=\Pi_{m=1}^k \lt(\rho\lt(\lt|z-w_m\rt|\rt)\frac{z-w_m}{\lt|z-w_m\rt|}\rt)^k.
\end{equation}
If $\rho(x)=x$ this is just a holomorphic function with order $k$ zero at $\lt\{w_1,w_2,\dots w_k\rt\}$. The idea is to 
create a function whose gradient close to an annulus of radius $1$ is very small. And whose gradient 
in the inside of the annulus and the outside of the annulus is large. 

Specifically we want estimates of the form 
\begin{equation}
\label{exeqq1}
\int_{B_{1-h}} \lt|Dh\rt| dz= O(1)\text{ and }\int_{B_2(0)\backslash B_{1+h}} \lt|Dh\rt| dz= O(1).
\end{equation}  
And for $\ep<<h$
\begin{equation}
\label{exeqq2}
\int_{B_{1+h}\backslash B_{1-h}} \lt|Dh\rt| dz\leq \ep.
\end{equation}
Now defining 
\begin{equation}
w(z):=\lt\{\begin{array}{ccc} h(z)-\Xint{-}_{\partial B_{1-h}} h\; dH^1 x & \text{ for } & z\in B_{1-h}\\
  l_R\circ h(z)-\Xint{-}_{\partial B_{1+h}} l_R \circ h\; dH^1 x & \text{ for } & z\in B_{2}\backslash B_{1+h}
\end{array}\rt.
\end{equation}
We can interpolate across $B_{1+h}\backslash B_{1-h}$ to create a function $\ti{w}$ with the property 
that 
\begin{equation}
D\ti{w}(z):=\lt\{\begin{array}{ccc} Dh(z) & \text{ for } & z\in B_{1-h}\\
  R Dh(z) & \text{ for } & z\in B_{2}
\end{array}\rt.
\end{equation}
and $\|\D\ti{w}\|_{L^{\infty}(B_{1+h}\backslash B_{1-h})}\leq c\ep$. If $h$ could be showed to 
be Quasiregular then we can use the method of \cite{astfar} "project" $\ti{w}$ onto the space of Quasiregular mappings to obtain a Quasiregular 
mappings with the properties required. So the main obstacle is to obtain a Quasiregular mapping that 
has properties (\ref{exeqq1}), (\ref{exeqq2}).  

Let
\begin{equation}
\label{ooeqq1}
G(z):=\Pi_{m=1}^k \lt(\rho(\lt|z-w_m\rt|)\rt)^k=e^{\frac{k^2}{2\pi}\lt(\sum_{m=1}^k \frac{2\pi}{k} \log\lt(\rho\lt(\lt|z-w_m\rt|\rt)\rt)\rt)}.
\end{equation}
Take $z=1$. Then 
\begin{eqnarray}
\lt|z-w_m\rt|&=&\sqrt{\lt(\lt(1-\cos\lt(\frac{2\pi m}{k}\rt)\rt)^2+\lt(\sin\lt(\frac{2\pi m}{k}\rt)\rt)^2\rt)}\nn\\
&=&\sqrt{2\lt(1-\cos\lt(\frac{2\pi m}{k}\rt)\rt)}.
\end{eqnarray}
So 
\begin{eqnarray}
\sum_{m=1}^k \frac{2\pi}k{ \log\lt(\rho\lt(\lt|1-w_m\rt|\rt)\rt)}&=&\sum_{m=1}^k 
\frac{2\pi}{k} \log\lt(\rho\lt(\sqrt{2\lt(1-\cos\lt(\frac{2\pi m}{k}\rt)\rt)}\rt)\rt)\nn\\
&\rightarrow& \int_{0}^{2\pi} \log\lt(\rho\lt(\sqrt{2\lt(1-\cos\lt(x\rt)\rt)}\rt)\rt) dx\nn\\
&=& \int_{0}^2 \log\lt(\rho\lt(r\rt)\rt)\frac{4}{\sqrt{4-r^2}} dr\nn\\
&=:&A_{\rho}.
\end{eqnarray}
Since $1$ is a typical point on $\partial B_1$ by symmetry of $z_1,z_2,\dots z_m$ so we have 
\begin{equation}
\label{exameq1}
\inf_{z\in \partial B_1(0)} G(z)\leq c e^{\frac{k^2}{2\pi} A_{\rho}}. 
\end{equation}
Let $\varpi(x)=\sum_{m=1}^{k} k\log\lt(\rho\lt(\lt|z-w_m\rt|\rt)\rt)$. So 
\begin{eqnarray}
\label{eqe401}
\int_{B_1} \lt|G(z)\rt| dz&=&\int_{B_1} e^{\log\lt(\lt|G(z)\rt|\rt)} dz\nn\\
&\overset{(\ref{ooeqq1})}{=}&\int_{B_1} e\circ \varpi (z) dz
\end{eqnarray}
Since $e^x$ is convex by Jensen's inequality we know 
\begin{equation}
\label{eqex402}
e^{\lt(\int_{B_1} \varpi(z) dz\rt)}\leq \int_{B_1} e\circ \varpi(z) dz.
\end{equation}
Let 
\begin{equation}
\label{ooeqq60}
 B_{\rho}:=\int_{0}^2 2 r\cos^{-1}\lt(\frac{r}{2}\rt)\log(\rho(r)) dr. 
\end{equation}
And note 
\begin{eqnarray}
\label{ooeqq61}
\int_{B_1} \varpi(z) dz&=&\sum_{m=1}^k k\int_{B_1} \log\lt(\rho\lt(\lt|z-w_m\rt|\rt)\rt) dz\nn\\
&=&k^2 \int_{B_1}  \log\lt(\rho\lt(\lt|z-(-1,0)\rt|\rt)\rt) dz\nn\\
&=&k^2 \int_{0}^2 2 r\cos^{-1}\lt(\frac{r}{2}\rt)\log(\rho(r)) dr \nn\\
&\overset{(\ref{ooeqq60})}{=}& k^2 B_{\rho}.
\end{eqnarray}
$$
\int_{B_1} \lt|G(z)\rt| dz\overset{(\ref{ooeqq60}),(\ref{eqex402}),(\ref{eqe401})}{\geq} e^{k^2 B_{\rho}}
$$
Thus a counter example can be constructed by finding an increasing function $\rho$ that satisfies the following 
two inequalities 
\begin{equation}
\label{ooeqq65}
A_{\rho}=\int_{0}^2 2 r\cos^{-1}\lt(\frac{r}{2}\rt)\log(\rho(r)) dr>0\text{ and } B_{\rho}=\int_{0}^2 \log\lt(\rho\lt(r\rt)\rt)\frac{4}{\sqrt{4-r^2}} dr<0
\end{equation}
and for which function $G$ defined (\ref{ooeqq1}) forms a quasiregular mapping. These things will 
be addressed in forthcoming preprint \cite{lor39}.

\section{Appendix}

We will prove an estimate from \cite{fmul} where we track the constants explicitly. All the arguments are from 
\cite{fmul}. 

\begin{a5} 
\label{PP1}
Suppose $u\in W^{1,2}(B_1:\R^2)$ with $\int_{B_1} \mathrm{dist}^2(Du,SO(2)) dz\leq 1$ 
then there exists $R\in SO(2)$ such that 
\begin{eqnarray}
\label{breq90}
&~&\int_{B_{\frac{1}{4}}} \lt|Du-R\rt|^2 dz\nn\\
&~&\qd\qd \leq 5 \lt(\int_{B_1} \mathrm{dist}^2(Du,SO(2)) dx\rt)^{\frac{1}{4}}+2\lt(\int_{B_1} \mathrm{dist}^2(Du,SO(2)) dx\rt)^{\frac{1}{4}}\|Du\|^{\frac{1}{2}}_{L^2(B_1)}\nn\\
\end{eqnarray}
\end{a5}

\em Step 1. \rm We will show 
\begin{equation}
\label{beq1}
\lt|\mathrm{cof}(M)-M\rt|\leq 2\mathrm{dist}(M,SO(2))\text{ for any }M\in M^{2\times 2}.
\end{equation}

\em Proof of Step 1. \rm Let $R_M\in SO(2)$ be such that $\lt|M-R_M\rt|=\mathrm{dist}(M,SO(2))$. Note 
$\lt|\mathrm{cof}(M)-R_M\rt|=\mathrm{dist}(M,SO(2))$. So 
$\lt|\mathrm{cof}(M)-M\rt|\leq \lt|\mathrm{cof}(M)-R_M\rt|+\lt|R_M-M\rt|= 2\mathrm{dist}(M,SO(2))$. Which establishes (\ref{beq1}). \nl

\em Step 2. \rm For any $w\in W^{1,2}(B_1,\R^2)$ we will show 
\begin{equation}
\label{beq2}
\int_{B_1} \lt|\lt|Dw\rt|^2-2\rt| dx\leq \lt(\int_{B_1} \mathrm{dist}^2\lt(Dw,SO(2)\rt) dx\rt)^{\frac{1}{2}}\lt(\|Dw\|_{L^2(B_1)}+\sqrt{2\pi}\rt).
\end{equation}

\em Proof of Step 2. \rm For any $x\in B_1$ let $R_x\in SO(2)$ be such that $\lt|Dw(x)-R_x\rt|=\mathrm{dist}(Dw(x),SO(2))$. So 
\begin{eqnarray}
\int_{B_1} \lt|\lt|Dw(x)\rt|^2-2\rt| dx&=&\int_{B_1} \lt|\lt( \lt|Dw(x)\rt|-\lt|R_x\rt|\rt)\lt(\lt|Dw(x)\rt|+\sqrt{2}\rt)\rt| dx\nn\\
&\leq& \lt( \int_{B_1}\lt|Dw(x)-R_x\rt|^2 dx\rt)^{\frac{1}{2}}\lt( \int_{B_1}(\lt|Dw(x)\rt|+\sqrt{2})^2 dx\rt)^{\frac{1}{2}}\nn\\
&\leq& \lt( \int_{B_1}\mathrm{dist}^2(Dw(x),SO(n)) dx\rt)^{\frac{1}{2}}\lt(\|Dw\|_{L^2(B_1)}+\sqrt{2\pi}\rt)
\end{eqnarray}
which establishes (\ref{beq2}).\nl

\em Proof of Proposition completed. \rm Let $z:B_1\rightarrow \R^2$ be the solution of 
$$
\triangle z=\mathrm{div}\lt(\mathrm{cof}(Du)-Du\rt), z=0\text{ on }\partial B_1.
$$
So testing the equation with $z$ itself we have 
\begin{eqnarray}
\int_{B_1} \lt|D z\rt|^2 dx&=&\int_{B_1} \lt(\mathrm{cof}(Du)-Du\rt):Dz dx\nn\\
&\leq& \lt(\int_{B_1} \lt|\mathrm{cof}(Du)-Du\rt|^2 dx\rt)^{\frac{1}{2}}\lt(\int_{B_1} \lt|Dz\rt|^2 dx\rt)^{\frac{1}{2}}\nn\\
&\overset{(\ref{beq1})}{\leq}&2 \lt(\int_{B_1} \mathrm{dist}^2\lt(Du,SO(2)\rt) dx\rt)^{\frac{1}{2}}\|Dz\|_{L^2(B_1)}.\nn
\end{eqnarray}
So 
\begin{equation}
\label{breq14}
\int_{B_1} \lt|Dz\rt|^2 dx\leq 4\int_{B_1} \mathrm{dist}^2\lt(Du,SO(2)\rt) dx.
\end{equation}
Let 
\begin{equation}
\label{ffefq41}
w=u-z. 
\end{equation}
Now using the identity 
$$
\frac{1}{2} \triangle (\lt|\nabla f\rt|^2)=\nabla f\cdot \triangle \nabla f+\lt|\nabla^2 f\rt|^2\text{ for any scalar valued function}f\in C^2. 
$$
So as $w$ is a vector valued function both of whose co-ordinates are harmonic we have 
\begin{equation}
\label{breq95}
\frac{1}{2} \triangle (\lt|D w\rt|^2-2)=\lt|D^2 w\rt|^2.
\end{equation}
Let $\eta\in C_0(B_1)$ be such that $\eta=1$ on $B_{\frac{1}{2}}$ and $\|D^2 \eta\|_{L^{\infty}(B_1)}\leq 8$. So 
\begin{eqnarray}
\int_{B_1} \lt|D^2 w\rt|^2 \eta dz&\overset{(\ref{breq95})}{=}&\int_{B_1} \frac{1}{2} \triangle \lt(\lt|D w\rt|^2-2\rt)\eta dx\nn\\
&=&\int_{B_1} \frac{1}{2}\lt(\lt|Dw\rt|-2\rt)\triangle \eta dx\nn\\
&\leq&\frac{1}{2}\sup_{B_1} \lt|\triangle \eta\rt|\int_{B_1} \lt|\lt|Dw\rt|^2-2\rt| dx\nn\\
&\leq&4\int_{B_1} \lt|\lt|Du\rt|^2-2Du:Dz+\lt|Dz\rt|^2-2\rt| dx\nn\\
&\leq& 4\lt(\int_{B_1} \lt|\lt|Du\rt|^2-2\rt|dx+  \int_{B_1} \lt|Dz\rt|^2 dx+2\lt(\int_{B_1} \lt|Dz\rt|^2 dx\rt)^{\frac{1}{2}}  \lt(\int_{B_1} \lt|Du\rt|^2 dx\rt)^{\frac{1}{2}} \rt)\nn\\
&\overset{(\ref{beq2}),(\ref{breq14})}{\leq}& 4\lt(\int_{B_1} \mathrm{dist}^2\lt(Du(x),SO(2)\rt) dx\rt)^{\frac{1}{2}}\lt(\|Du\|_{L^2(B_1)}+\sqrt{2\pi}\rt)\nn\\
&~&\qd+16\int_{B_1} \mathrm{dist}^2\lt(Du(x),SO(2)\rt) dx\nn\\
&~&\qd+4\lt(\int_{B_1} \mathrm{dist}^2\lt(Du(x),SO(2)\rt) dx\rt)^{\frac{1}{2}}\|Du\|_{L^2(B_1)}.\nn\\
&\leq& 8\lt(\int_{B_1} \mathrm{dist}^2\lt(Du(x),SO(2)\rt) dx\rt)^{\frac{1}{2}}\|Du\|_{L^2(B_1)}\nn\\
&~&\qd +27\lt(\int_{B_1} \mathrm{dist}^2\lt(Du(x),SO(2)\rt) dx\rt)^{\frac{1}{2}}.
\end{eqnarray}
So 
\begin{eqnarray}
\label{breq96}
\lt(\int_{B_{\frac{1}{2}}} \lt|D^2 w\rt|^2 dx\rt)^{\frac{1}{2}}&\leq& 2\sqrt{2}\lt(\int_{B_1} \mathrm{dist}^2\lt(Du(x),SO(2)\rt) dx\rt)^{\frac{1}{4}}\|Du\|_{L^2(B_1)}^{\frac{1}{2}}\nn\\
&~&\qd+3\sqrt{3}\lt(\int_{B_1} \mathrm{dist}^2\lt(Du(x),SO(2)\rt) dx\rt)^{\frac{1}{4}}.
\end{eqnarray}
Note $\int_{B_{\frac{1}{2}}} \lt|D^2 w\rt| dx\leq \lt( \int_{B_{\frac{1}{2}}} \lt|D^2 w\rt|^2 dx\rt)^{\frac{1}{2}}\sqrt{\frac{\pi}{4}}$. Let 
$y\in B_{\frac{1}{4}}$, by the mean value theorem 
$$
D^2 w(y)=\Xint{-}_{B_{\frac{1}{4}}(y)} D^2 w(x) dx.
$$
So 
\begin{eqnarray}
\lt|D^2 w(y)\rt|&\leq& \lt(\frac{\pi}{16}\rt)^{-1}\int_{B_{\frac{1}{4}}(y)} \lt|D^2 w(x)\rt| dx\nn\\
&\overset{(\ref{breq96})}{\leq}& \frac{16}{\pi}\times 2\sqrt{2}\lt(\int_{B_1} \mathrm{dist}^2\lt(Du(x),SO(2)\rt) dx\rt)^{\frac{1}{4}}\|Du\|_{L^2(B_1)}^{\frac{1}{2}}\nn\\
&~&\qd\qd+\frac{16}{\pi}\times 3\sqrt{3}\lt(\int_{B_1} \mathrm{dist}^2\lt(Du(x),SO(2)\rt) dx\rt)^{\frac{1}{4}}.
\end{eqnarray}
So 
\begin{eqnarray}
\label{fggf1}
\|D^2 w(y)\|_{L^{\infty}(B_{\frac{1}{4}})}&\leq& \frac{32\sqrt{2}}{\pi}\lt(\int_{B_1}\mathrm{dist}^2\lt(Du(x),SO(2)\rt) dx\rt)^{\frac{1}{4}}\|Du\|_{L^2(B_1)}^{\frac{1}{2}}\nn\\
&~&\qd\qd\qd+\frac{48\sqrt{3}}{\pi}\lt(\int_{B_1} \mathrm{dist}^2\lt(Du(x),SO(2)\rt) dx\rt)^{\frac{1}{4}}.
\end{eqnarray}
Let $x_0\in B_{\frac{1}{4}}$. Thus 
\begin{eqnarray}
\label{breq101}
\sup\lt\{\lt|Dw(x)-Dw(x_0)\rt|:x\in B_{\frac{1}{4}}\rt\}&\leq& \frac{16\sqrt{2}}{\pi}\lt(\int_{B_1} \mathrm{dist}^2\lt(Du(x),SO(2)\rt) dx\rt)^{\frac{1}{4}}\|Du\|_{L^2(B_1)}^{\frac{1}{2}}\nn\\
&~&+\frac{24\sqrt{3}}{\pi}\lt(\int_{B_1} \mathrm{dist}^2\lt(Du(x),SO(2)\rt) dx\rt)^{\frac{1}{2}}.
\end{eqnarray}
Now
\begin{equation}
\label{ffeqf101}
\int_{B_{\frac{1}{4}}} \lt|Dz \rt| dz\leq \lt(\int_{B_{\frac{1}{4}}} \lt|Dz \rt|^2 dz\rt)^{\frac{1}{2}}\frac{\sqrt{\pi}}{4}
\overset{(\ref{breq14})}{\leq}\frac{\sqrt{\pi}}{2}\lt(\int_{B_1} \mathrm{dist}^2(Du,SO(2)) dz\rt)^{\frac{1}{2}}.
\end{equation}
And
\begin{eqnarray}
\label{breq102}
\int_{B_{\frac{1}{4}}} \lt|\na u(x)-\na w(x_0)\rt| dx&\leq& \int_{B_{\frac{1}{4}}} \lt|\na u(x)-\na w(x)\rt| dx
+\int_{B_{\frac{1}{4}}} \lt|\na w(x)-\na w(x_0)\rt| dx\nn\\
&\overset{(\ref{ffefq41})}{\leq}& \int_{B_{\frac{1}{4}}} \lt|\na z(x)\rt| dx+\frac{\pi}{16}\|D w-D w(x_0)\|_{L^{\infty}(B_{\frac{1}{4}})}\nn\\
&\overset{(\ref{breq101})}{\leq}&\int_{B_{\frac{1}{4}}} \lt|D z(x)\rt| dx+\sqrt{2}\lt(\int_{B_1} \mathrm{dist}^2\lt(Du(x),SO(2)\rt) dx\rt)^{\frac{1}{4}}\|Du\|_{L^2(B_1)}^{\frac{1}{2}}\nn\\
&~&\qd+\frac{3\sqrt{3}}{2}\lt(\int_{B_1}\mathrm{dist}^2\lt(Du(x),SO(2)\rt) dx\rt)^{\frac{1}{2}}\nn\\
&\overset{(\ref{ffeqf101})}{\leq}&\lt(\frac{3\sqrt{3}}{2}+\frac{\sqrt{\pi}}{2}\rt) \lt(\int_{B_1} \mathrm{dist}^2\lt(Du(x),SO(2)\rt) dx\rt)^{\frac{1}{4}}\nn\\
&~&\qd+\sqrt{2}\lt(\int_{B_1} \mathrm{dist}^2\lt(Du(x),SO(2)\rt) dx\rt)^{\frac{1}{4}}\|Du\|_{L^2(B_1)}^{\frac{1}{2}}.
\end{eqnarray}

Recall $w=u-z$. So 
\begin{eqnarray}
\lt(\int_{B_1} \mathrm{dist}^2\lt(Dw,SO(2)\rt) dx\rt)^{\frac{1}{2}}&\leq& 
\lt(\int_{B_1} \mathrm{dist}^2\lt(Du,SO(2)\rt) dx\rt)^{\frac{1}{2}}+\|Dz\|_{L^2(B_1)}\nn\\
&\overset{(\ref{breq14})}{\leq}&3\lt(\int_{B_1} \mathrm{dist}^2\lt(Du,SO(2)\rt) dx\rt)^{\frac{1}{2}}.\nn
\end{eqnarray}
Hence $\int_{B_1} \mathrm{dist}(Dw,SO(2)) dx\leq 3\sqrt{\pi}\lt(\int_{B_1} \mathrm{dist}^2(Du,SO(2)) dx\rt)^{\frac{1}{2}}$. So there 
must exist $x_0\in B_1$ such that 
\begin{equation}
\label{breq99}
\mathrm{dist}(Dw(x_0),SO(2))\leq \frac{3}{\sqrt{\pi}}\lt(\int_{B_1} \mathrm{dist}^2(Du,SO(2)) dx\rt)^{\frac{1}{2}}.
\end{equation}
Let $R\in SO(2)$ be such that $\lt|Dw(x_0)-R\rt|=\mathrm{dist}(Dw(x_0),SO(2))$. By (\ref{breq102}), (\ref{breq99}) we have that 
\begin{eqnarray}
\label{breq300}
\int_{B_{\frac{1}{4}}} \lt|Du(x)-R\rt| dx&\leq&5\lt(\int_{B_1} \mathrm{dist}^2(Du,SO(2)) dx\rt)^{\frac{1}{4}}\nn\\
&~&\qd\qd+2\lt(\int_{B_1} \mathrm{dist}^2(Du,SO(2)) dx\rt)^{\frac{1}{4}}\|Du\|^{\frac{1}{2}}_{L^2(B_1)}.
\end{eqnarray}

\end{document}